\documentclass{article}

\usepackage{graphicx} 
\usepackage{todonotes}
\usepackage{mathrsfs,bm}
\usepackage[a4paper, total={6in, 9in}]{geometry}
\usepackage{mathtools}
\usepackage{enumitem}
\usepackage{amsmath}
\usepackage{amsthm}
\usepackage{amssymb}
\usepackage{graphicx}
\usepackage[english]{babel}
\usepackage[dvipsnames]{xcolor}
\usepackage[autostyle=true]{csquotes}
\usepackage{tikz-cd}

\newtheorem{theorem}{Theorem}[subsection]
\newtheorem{corollary}[theorem]{Corollary}
\newtheorem{lemma}[theorem]{Lemma}
\newtheorem{proposition}[theorem]{Proposition}

\theoremstyle{definition}
\newtheorem{remark}[theorem]{Remark}

\usepackage[dvipsnames]{xcolor}
\usepackage[colorlinks=true]{hyperref}
\hypersetup{urlcolor=RubineRed,linkcolor=RoyalBlue,citecolor=ForestGreen}

\usepackage[
backend=biber,
style=alphabetic,
sorting=nyt
]{biblatex}
\renewbibmacro{in:}{}

\usepackage{sectsty}
\sectionfont{\centering}
\numberwithin{equation}{section}

\title{Non-archimedean generalisations of the Kuznetsov formulas}
\author{Matteo Di Scipio}
\date{November 2026}

\begin{document}

\maketitle
\abstract{\noindent In this paper we develop non-archimedean analogues of the Petersson-Kuznetsov formulas. Specifically, choosing a particular test function on either the geometric or spectral side of the formula, one may express the other side in terms of integral transforms. Such transforms involve representation theoretic Bessel functions attached to irreducible unitary admissible representations of $\text{PGL}_2(\mathbb{Q}_p)$.}

\setcounter{tocdepth}{2}
\tableofcontents
\pagebreak

\section{Introduction}
\subsection{Motivation and main results}
    In the late 70s, Kuznetsov \cite{Kuznetsov1980} and Bruggeman introduced summation formulae relating sums of pairs of Fourier coefficients of Maass forms (the \emph{spectral} side) to sums of Kloosterman sums (the \emph{geometric} side), which have since become part of the toolkit of analytic number theorists to study several questions concerning automorphic forms.
In a representation theoretic framework, given a well behaved complex valued function $h_\infty$ on the unitary dual of $G(\mathbb{R})$, where $G=\operatorname{PGL}_2$, a positive integer $N$ and integers $mn>0$ one has
\begin{align}\label{eq:forward-kuznetsov}
    \sum_{\pi \ \text{cuspidal}}h_\infty(\pi) & \sum_{\varphi \in \mathcal{B}(\pi^{K_0(N)})}a_{u_\varphi}(m) \overline{a_{u_\varphi}(n)} \ + \ (\text{continuous term}) \notag \\
    & =\delta_{m=n}\int_{G(\mathbb{R})^\wedge}h_\infty(\pi)d\widehat{\mu}(\pi)+\sum_{c\in N\mathbb{N}}\frac{S(m,n,c)}{c}H_\infty\left(\frac{4\pi\sqrt{mn}}{c}\right)
\end{align}
where the outer sum on the LHS runs over all cuspidal automorphic representations of $G$ over $\mathbb{Q}$, $\mathcal{B}(\pi^{K_0(N)})$ is an orthonormal basis for the right-$K_0(N)$-invariant subspace of $\pi$, $a_{u_\varphi}(m)$ is the Fourier coefficient of a modular/Maass form underlying $\varphi$ as described in Section \ref{section:inner-products} and $H_\infty$ is an integral transform of $h_\infty$ given by
\begin{equation}\label{def:archim-H-transform}
    H_\infty(x)=\int_{G(\mathbb{R})^\wedge}h_\infty(\pi)J_{\pi}(x)d\widehat{\mu}(\pi).
\end{equation}
Here $J_\pi$ are the kernels
\begin{equation}\label{def:J_pi}
    J_\pi(x)=\begin{cases}
        2\pi i^{-\kappa}J_{\kappa-1}(x) & \pi\simeq \mathcal{D}(\kappa), \\
        \pi i\frac{J_{2it}(x)-J_{-2it}(x)}{\sinh(\pi t)} & \pi \simeq \pi(\text{sgn}^\epsilon |\cdot|^{it},\text{sgn}^\epsilon |\cdot|^{-it}),
    \end{cases}
\end{equation}
and $d\widehat{\mu}(\pi)$ is the Plancherel measure on the unitary dual as described in Section \ref{section:basics-on-local-theory}.

In the same paper Kuznetsov also introduced a companion formula where one chooses a test function appearing on the geometric side and the integral transform appears on the spectral side instead. For $f_\infty \in C_c^\infty(\mathbb{R}_+)$ and $m,n,N$ as before, one has 
\begin{align}\label{eq:backward-kuznetsov}
    \sum_{c\in N\mathbb{N}}\frac{S(m,n,c)}{c} & f_\infty\left(\frac{4\pi\sqrt{mn}}{c}\right) \notag
    \\ & =\sum_{\pi \ \text{cuspidal}}\widehat{H}_{\infty}(\pi)\sum_{\varphi \in \mathcal{B}(\pi^{K_0(N)})}a_{u_\varphi}(m)\overline{a_{u_\varphi}(n)}+(\text{continuous term})
\end{align}
where
\begin{equation}\label{def:archim-backward-transform}
    \widehat{H}_\infty(\pi)=\frac{2}{\pi}\int_{0}^{\infty}J_\pi(x)f_\infty(x)d^\times x.
\end{equation}
and all other notation as before.

Although as discussed in \cite{Motohashi}, formulae \eqref{eq:forward-kuznetsov} and \eqref{eq:backward-kuznetsov} are essentially equivalent, the two versions have found plenty of applications in analytic number theory; for the former see for instance \cite{Conrey2000}, \cite[\S 26]{IwaniecKowalski2004} and \cite{Petrow2019} for applications to moments and non-vanishing, and for the latter see \cite{DeshouillersIwaniec1982} for bilinear forms of Kloosterman sums and improvements on the Selberg conjecture, \cite{Young2011} and \cite{Petrow2023} for moments of Dirichlet $L$-functions . Here and in the sequel, we refer to \eqref{eq:forward-kuznetsov} as the \emph{spectral-geometric} formula and  \eqref{eq:backward-kuznetsov} as the \emph{geometric-spectral} formula. Given the extensive literature on Bessel functions, the transforms $H_\infty$ and $\widehat{H}_\infty$ can be fairly well understood through standard methods such as integration by parts and stationary phase, allowing for strong bounds on sums of Fourier coefficients or Kloosterman sums. 

\bigbreak 
The goal of this paper is establishing generalisations of the Kuznetsov formulae that allows to select a $p$-adic test function on one side and express the other side in terms of integral transforms. For the spectral-geometric formula, given a suitably nice complex valued function $h_p$ defined on the dual $G(\mathbb{Q}_p)^\wedge$, and used to pick out certain families of representations, we seek a kernel $k(\pi,x)$ so that we have an integral transform $H_p$ appearing on the geometric side defined by
\begin{equation}
    h_p \quad \longmapsto  \quad H_p(x):=\int_{G(\mathbb{Q}_p)^\wedge}k(\pi,x) h_p(\pi)d\widehat{\mu}(\pi)
\end{equation}
and similarly, for the geometric-spectral formula, given a nice function $f \in C(\mathbb{Q}_p^\times)$ acting as a non-archimedean weight to the Kloosterman sums, we want a transform
\begin{equation}
    f \quad \longmapsto  \widehat{H}_p(\pi):=\int_{\mathbb{Q}_p^\times}k(\pi,x)f(x)d^\times x
\end{equation}
We will show the kernels $k(\pi,x)$ arise from representation theory, and are related to the action of the Weyl element in the Kirillov model.

\subsubsection{The space of test functions and the transforms} 
We invite the reader to go through Section \ref{section:basics-on-local-theory} and the beginning of Section \ref{section:non-arch spectral decomp} for the exact representation theoretic notions. \\
\noindent Spectral-geometric direction: for $F=\mathbb{R}$ or $\mathbb{Q}_p$ the unitary dual $G(F)^\wedge$ is equipped with the Fell topology and Plancherel measure $\widehat\mu$. Denote by $G(F)^\wedge_{\operatorname{gen}}$ the topological subspace of generic representations. Let $\mathcal{S}(G(F)^\wedge)$ the space of complex valued functions supported on the generic irreducible representations satisfying the following conditions: \\
\noindent If $F=\mathbb{R}$, then $h_\infty \in \mathcal{S}(G(\mathbb{R})^\wedge)$ is such that, for some $\delta>0$,
\begin{itemize}[itemsep=0em]
    \item $h_\infty(\mathcal{D}(\kappa))\ll \kappa^{-2-\delta}$ on the discrete series,
    \item On principal series $h_\infty(\pi(\text{sgn}^\epsilon |\cdot|^{it},\text{sgn}^\epsilon |\cdot|^{-it}))=h(t)$ where $h(t)$ satisfies
    \begin{equation}
    \begin{cases}
        h(t) \text{ is even}, \\
        h(t) \text{ is holomorphic on the strip } |\text{Im}(t)|\leq 13,\\
        h(t)\ll (1+|t|)^{-2-\delta}.
    \end{cases}
\end{equation}
\end{itemize}
If instead $F=\mathbb{Q}_p$, and $\pi$ an irreducible unitary generic representation of $\operatorname{PGL}_2(\mathbb{Q}_p)$ denote by $[\pi]\subseteq \operatorname{PGL}_2(\mathbb{Q}_p)^\wedge_{\operatorname{gen}}$ the connected component in the unitary generic dual containing $\pi$. Note that if $\pi$ is a supercuspidal then $[\pi]$ is a singleton.
Take $\mathcal{S}(\operatorname{PGL}_2(\mathbb{Q}_p)^\wedge)$ to be the space of functions on $\operatorname{PGL}_2(\mathbb{Q}_p)^\wedge_{\operatorname{gen}}$ spanned by
\begin{itemize}[itemsep=0em]
    \item The indicator functions $1([\pi])$ for any $\pi$ that is either supercuspidal or an induced representation $\pi(\chi,\chi^{-1})$ where $\chi|_{\mathbb{Z}_p^\times}$ is not quadratic .
    \item The function $h$ supported on spherical representations such that
    $$h(\pi(|\cdot|^{i\theta/\log p},|\cdot|^{-i\theta/\log p}))=P(\cos \theta), \quad \text{for some } P\in \mathbb{C}[X]$$
\end{itemize}
The extra flexibility on the unramified representations is due to Satake isomorphism, and might be of interest to establish Sato-Tate type equidistribution laws, such as in \cite[\S 10]{KL-kuznetsov}. The requirement $\chi^2|_{\mathbb{Z}_p}\neq 1$ is present for technical reasons, and it is essentially due to the non-existence ofcompactly supported newform projectors on the $c(\pi)=1$ and $2$ lines. We refer the reader to Proposition \ref{prop:proj-nonexistence}.\\
For $h_\infty \in \mathcal{S}(G(\mathbb{R})^\wedge)$ we then consider the transform $H_\infty$ as described in \eqref{def:archim-H-transform}, whereas if $h_p\in \mathcal{S}(\operatorname{PGL}_2(\mathbb{Q}_p)^\wedge)$ we define
\begin{equation}\label{def:non-archim-H_transform}
    H_p(x)=\frac{1}{\zeta_p(1)}\int_{G(\mathbb{Q}_p)^\wedge}j_{\pi,\psi_p}(x)\frac{h_p(\pi)}{\mathcal{L}_\pi(1)}d\widehat\mu(\pi)
\end{equation}
where $j_{\pi,\psi_p}$ is the $j$-Bessel function attached to $\pi$ and the additive character $\psi_p$ as defined in Section \ref{section:Bessel-functions}, and $\mathcal{L}_\pi(1)$ is a constant depending on the isomorphism class of $\pi$ as in \eqref{def:euler-factors}. \bigbreak

\noindent Geometric-spectral direction: we let $\mathcal{S}(F^\times):=C_c^\infty(F^\times)$ the space of smooth, compactly supported complex valued functions on $F^\times$. For $f_\infty \in \mathcal{S}(\mathbb{R}_+)$ we take $\widehat{H}_\infty$ as in \eqref{def:archim-backward-transform}. For $f \in \mathcal{S}(\mathbb{Q}_p^\times)$, further define the transform $\widehat{H}_p$ by
\begin{equation}\label{def:non-archim-backward-transform}
    \widehat H_p(\pi)=\mathcal{L}_\pi(1)\frac{\zeta_p(2)}{\zeta_p(1)}\int_{\mathbb{Q}_p^\times}f(u)j_{\pi,\psi_p}(u) d^\times u
\end{equation}
where we keep the same notation as above.

\begin{remark}
    To push the analogy between $H_\infty$ and $H_p$ further (as well as $\widehat{H}_\infty$ and $\widehat{H}_p$), it is worth pointing out that the kernels $J_\pi$ given in \eqref{def:J_pi} arise from archimedean representation theoretic Bessel functions as well. See the formulas in Section \ref{section:Bessel-functions}.
\end{remark}

\subsubsection{Statement of the results}
\begin{theorem}[Spectral-geometric formula]\label{thm:forward-formula} Fix $m,n \in \mathbb{Z}$ such that $mn>0$ and $p$ a prime such that $(mn,p)=1$. Let $h_\infty \in \mathcal{S}(\operatorname{PGL}_2(\mathbb{R})^\wedge)$  and $h_p \in \mathcal{S}(\operatorname{PGL}_2({\mathbb{Q}_p)}^\wedge)$. We have
\begin{align}
    \sum_{\substack{\pi \ \text{cuspidal}\\ C(\pi)|p^{\infty}}}h_\infty(\pi_\infty)h_p(\pi_p) & a_{u_\varphi}(m)\overline{a_{u_\varphi}(n)}+(\text{continuous term})\notag \\
    & = \delta_{m=n}\delta_\infty\delta_p + \sum_{\mu\in \mathbb{Q}_+} \frac{H(m,n,\mu) }{\mu^{1/2}}H_\infty\left(4\pi\sqrt{\frac{mn}{\mu}}\right)H_p\left(\frac{mn}{\mu}\right)
\end{align}
with both sides converging absolutely. Here
\begin{itemize}[itemsep=-0.2em]
    \item $C(\pi):=\prod_{p}p^{c(\pi_p)}$ is the conductor of the representation $\pi$,
    \item $a_{u_{\varphi}}(m)$ denotes the Fourier coefficient of the modular/Maass form underlying the $L^2$-normalised newform $\varphi \in \pi$, as in Section \ref{section:inner-products},
    \item the continuous term is given in \eqref{def:foward-thm-cts-term}, which vanishes if $h_p$ is supported only on supercuspidal representations or $h_\infty$ is supported only on discrete series,
    \item $H_\infty$ and $H_p$ denote the transforms given in \eqref{def:archim-H-transform} and \eqref{def:non-archim-H_transform}, and $H(m,n,\mu)$ are the Kloosterman sums
    $$H(m,n,\mu)=\begin{cases}
        S(\overline{p^k}m,n,c) & \text{if } \mu=c^2p^k, (c,p)=1, \overline{p}p\equiv 1 \text{ mod }c,\\ 0 & \text{otherwise},
    \end{cases}$$
    where $c \in \mathbb{N}$, and $k\in \mathbb{Z}$. The $\mu$-sum runs over a discrete subset of $\mathbb{Q}_+$.
    \item The archimedean and non-archimedean diagonal terms $\delta_\infty$ and $\delta_p$ are given by
    $$\delta_\infty=\int_{G(\mathbb{R})^\wedge}h_\infty(\pi)d\widehat{\mu}(\pi) \qquad \delta_p=\int_{G(\mathbb{Q}_p)^\wedge}\frac{h_p^*(\pi)}{\mathcal{L}_\pi(1)}d\widehat{\mu}(\pi)$$
    where $h^*(\pi)=h(\pi)$ if $c(\pi)\geq 1$ while for spherical representations $h^*_p$ is as in Lemma \ref{lemma:h-star}. 
\end{itemize}
\end{theorem}

\noindent To see how one may recover the classical Kuznetsov formula, see Remark \ref{remark:recover-classical-formula}.\\
We shall prove the Theorem \ref{thm:forward-formula} via the relative trace formula as in \cite{KL-kuznetsov} and \cite{KnigthlyLiPetersson}, by integrating an automorphic kernel attached to some test function $f_\infty\otimes f$ with $f_\infty\in C_c(G(\mathbb{R}))$ and $f\in C_c(G(\mathbb{A}_\text{fin}))$. It is interesting to note that the choice of archimedean/non-archimedean test function is suppressed in the final formula.

\noindent Similar formulas have previously appeared in the literature, see for instance \cite{Hu2024} and \cite{generalised}, but with some key differences. Firstly the geometric side is not expressed in terms of integral transforms $H_p$ of the test function $h_p$, but in terms of orbital integrals of a function $f_p \in C_c^\infty(G(\mathbb{Q}_p))$ underlying $h_p$ (see Remark \ref{remark:h-to-f} and Proposition \ref{prop:h_p-to-f_p}). In particular their results rely on being able to study $ f_p$, whereas very little information is required in our case. Secondly, given a supercuspidal representation $\sigma$ and denoting $\eta$ the unramifed character twist of $\mathbb{Q}_p^\times$, the authors study indicators on the neighborhood $[\sigma]_{\operatorname{HPY}}$ given by
$$[\sigma]_{\operatorname{HPY}}=\begin{cases}
    \{\sigma\} & p>2 \text{ and } c(\sigma) \text{ even}, \\
    \{ \sigma,\sigma\otimes \eta\} & p> 2  \text{ and } c(\sigma) \text{ odd},
\end{cases}$$
and consisting of two or three supercuspidal representations if $p=2$ (depending on the ramification of the underlying field extension $E/\mathbb{Q}_p$) . Compare this with Theorem \ref{thm:forward-formula} where we are able to pick out singletons independently of the conductor exponent and base prime. Our formula may be viewed as an analogue of the method of matrix pseudo-coefficients used in \cite{palm2012} used to pick out representations in the Selberg trace formula setting  to obtain hybrid Weyl laws (Thm 3.2.1 loc. cit.).
We may in fact refine the Theorem \ref{thm:forward-formula} to obtain a more explicit version:
\begin{corollary}\label{cor:supercuspidal-selection}
    Let $\sigma$ be a supercuspidal representation of $G(\mathbb{Q}_p)$ (such that $c(\pi)\geq 9$ if $p=2)$ with associated indicator function $h_\sigma\in \mathcal{S}(G(\mathbb{Q}_p)^\wedge)$. Let $h_\infty \in \mathcal{S}(G(\mathbb{R})^\wedge)$  and $m,n\in \mathbb{Z}$ such that $mn>0$ and $(mn,p)=1$, then
    \begin{align*}
        \sum_{\substack{\pi \ \text{cuspidal}\\ C(\pi)|p^\infty }}& h_\infty(\pi_\infty)  h_\sigma(\pi_p)a_{u_\varphi}(m)\overline{a_{u_\varphi}(n)}\\ &=\delta_{m=n}\delta_{\infty}\frac{d_\sigma}{(1-p^{-1})} +\sum_{\mu \in C(\sigma)\mathbb{Z}}\frac{S(\overline{p^{e_p(\mu)}}m,n,c)}{\mu^{1/2}}H_p\left(\frac{mn}{\mu}\right)H_\infty\left(4\pi \sqrt{\frac{mn}{\mu}}\right)
    \end{align*}
    with the following notation: 
    \begin{itemize}[itemsep=-0.2em]
        \item $C(\sigma):=p^{c(\sigma)}$ is the conductor of $\sigma$ and $d_\sigma$ is its formal degree, which if $\sigma$ is twist-minimal is given by
        $$d_\sigma=\begin{cases} (p-1)p^{r-1} & c(\sigma)=2r, \\
        \tfrac{p^2-1}{2}p^{r-1} & c(\sigma)=2r+1.    
        \end{cases}$$
        \item The transform $H_p$ is given by
        $ H_p(x)=d_\sigma j_{\sigma,\psi_p}(x) $ where the $j$-Bessel function is as in equation \eqref{eq:baruch-snitz},
        \item All other notation as before.
    \end{itemize}
\end{corollary}
\begin{remark} For the formulas for $d_\sigma$, see \cite[\S 5.11]{Carayol1984}. If $p$ is odd all supercuspidal representations of $G(\mathbb{Q}_p)$ are twist minimal, whereas if $p=2$ then the twist-minimal ones are those for which $c(\pi)=2$ or odd. If $(mn,p)\neq 1$ and $c(\pi)\geq 2$, then Atkin-Lenher theory implies that $a_{u_\varphi}(m)\overline{a_{u_\varphi}(n)}=0$. For a discussion on the geometric side in this case, we refer the reader to Lemma \ref{lemma:supercuspida-trivial-bound}. The restrictions on $\sigma$ when $p=2$ is so that only dihedral representations are considered, though it is worth remarking that Theorem \ref{thm:forward-formula} works for any choice of $\sigma$ and with a little more effort, refinements in those cases should be possible.  
\end{remark}

\noindent Note that if $\sigma$ is a supercuspidal representation of odd conductor exponent with unramified quadratic twist $\sigma\otimes \eta$, the two representations are not isomorphic and their $\varepsilon$-factors obey the relation $\varepsilon(1/2,\sigma,\psi)=-\varepsilon(1/2,\sigma\otimes \eta,\psi)$, so unlike \cite{generalised} Corollary \ref{cor:supercuspidal-selection} allows to select families of automorphic forms with fixed root number. This is particularly useful for the study of averages of (shifted) moments over families of automorphic forms, which is often carried out via the approximate functional equation of \cite[Thm 5.3]{IwaniecKowalski2004} which is dependent on the root number.

\noindent Another application of Corollary \ref{cor:supercuspidal-selection} is the non-vanishing of Fourier coefficients:
\begin{corollary}\label{cor:non-vanishing-coeff}
    Let $\sigma$ be a supercuspidal representation of conductor $N=p^r$, where $r$ is odd if $p=2$.  Let $\kappa\geq 2$ even and $m \in \mathbb{N}$ coprime to $p$. For $N$ sufficiently large there exist newforms $f \in S_\kappa^{\operatorname{new}}(N)$ such that $a_f(m)\neq 0$ and the underlying automorphic representation satisfies $(\pi_f)_p \simeq \sigma$.
\end{corollary}

We now turn to the geometric-spectral formula. In order to state the analogue of \eqref{eq:backward-kuznetsov} in its most general form, we introduce the following definitions. We say $f \in C_c^\infty(G(\mathbb{Q})_p)$ is a \emph{newform projector} if for all generic $\pi$ in the unitary dual, the operator $\pi(f)$ either projects on the newform line or it is identically zero. For a cofinite set $S$ of primes, we say the family of functions $\{f_p \in C_c^\infty(G(\mathbb{Q}_p))\}_{p\in S}$ is \emph{admissible} if:
\begin{enumerate}
    \item for all $p$, the function $f_p$ is either a newform projector or is of the form $\nu(p^c)1_{K_0(p^c)}$, where $\nu(n)=n\prod_{p|n}(1+p^{-1})$,
    \item $f_p=1_{K_p}$ for almost all primes in $S$.
\end{enumerate}
For such a family we construct the \emph{global generalised Kloosterman sum}
\begin{equation}
    H^S(m,n,\mu)=\prod_{p\in S}H_p(m,n,\mu), \qquad (m,n\in \mathbb{Z},\mu \in \mathbb{Q}_+)
\end{equation}
where $H_p(m,n,\mu)$ is the local generalised Kloosterman sum underlying $f_p$ given by
\begin{equation}\label{def:local-kloosterman-sum}
    H_p(m,n,\mu)=\iint_{\mathbb{Q}_p\times \mathbb{Q}_p}f_p(\begin{pmatrix}
        1 & -t_1 \\ & 1
    \end{pmatrix}\begin{pmatrix}
        & -\mu^{-1} \\ 1 &
    \end{pmatrix}\begin{pmatrix}
        1 & t_2 \\ & 1
    \end{pmatrix})\psi_p(mt_1-nt_2)dt_1dt_2
\end{equation}
Note that $H^S(m,n,\mu)$ in Theorem \ref{thm:forward-formula} is the Kloosterman sum for the family $\{1_{K_\ell}\}_{\ell\neq p}$.

\begin{remark}\label{remark:h-to-f}
    In fact, given $h_p \in \mathcal{S}(G(\mathbb{Q}_p)^\wedge)$, it is possible to express $H_p(x)$ as a local generalised Kloosterman sum for a test function $f_p \in C_c^\infty(G(\mathbb{Q}_p))$ such that the map
    $$h_p \longmapsto f_p \longmapsto (\pi \mapsto \operatorname{trace} \pi(f_p)) $$
    is the identity on $\mathcal{S}(G(\mathbb{Q}_p)^\wedge)$. This allows to recover the Kuznetsov-type formulas of \cite{Hu2024} and \cite{generalised}. See Propositions \ref{prop:h_p-to-f_p} and \ref{prop:kloosterman-to-transform}.
\end{remark}

\begin{theorem}[Geometric-spectral formula]\label{thm:backward-formula} Fix a finite prime $p$ and let $m,n\in \mathbb{Z}$ such that $mn>0$ and $(mn,p)=1$. Let $\{f_\ell\}_{\ell\neq p}$ be an admissible family of test functions. For $f_\infty \in \mathcal{S}(\mathbb{R}_+)$ and $f_p \in \mathcal{S}(\mathbb{Q}_p^\times)$, we have the spectral expansion 
    \begin{align*}
        \sum_{\mu\in \mathbb{Q}_+}\frac{H^S(m,n,\mu)}{\sqrt{\mu}|\mu|_p} & f_\infty\left(4\pi\sqrt{\frac{mn}{\mu}}\right)f_p\left(\frac{mn}{\mu}\right)\\
        & =\sum_{\pi \ \text{cuspidal}}\widehat{H}_\infty(\pi_\infty)\widehat{H}_p(\pi_p)\sum_{\varphi \in \mathcal{B}(\pi)^*}a_{u_\varphi}(m)\overline{a_{u_\varphi}(n)}+(\text{continuous term})
    \end{align*}
    Here:
    \begin{itemize}[itemsep=-0.2em]
        \item Given $\pi\simeq \bigotimes_v\pi_v $, the transforms $\widehat{H}_{\infty}(\pi_\infty)$ and $\widehat{H}_{p}(\pi_p)$ are given by \eqref{def:archim-backward-transform} and \eqref{def:non-archim-backward-transform},
        \item the continuous term is described in \eqref{eq:backward-formula-cts-term},
        \item $H^S(m,n,x)$ is the generalised Kloosterman sum for the family $\{f_\ell\}_{\ell \neq p}$,
        \item $\mathcal{B}(\pi)^*$ is an orthonormal basis consisting of pure tensors $\bigotimes _v \varphi_v$ such that $\varphi_\infty,\varphi_p$ are the newforms for $\pi_\infty$ and $\pi_p$, and for $\ell\neq p$ the local vector $\varphi_\ell \in \operatorname{Im}(\pi(f_\ell))$ .   
    \end{itemize}
\end{theorem}

\noindent It is interesting to compare Theorem \ref{thm:backward-formula} with the spectral expansions established by Blomer and Milicevic in \cite[Thm 4]{Blomer_2015}. The latter was employed by the authors to study cancellation in weighted sums of Kloosterman sums: for $q$ a prime power and  $f: (\mathbb{Z}/q\mathbb{Z})^\times \to \mathbb{C} $ they obtain
\begin{align*}
    \sum_{(c,p)=1}\frac{f_\infty(c/X)}{c^{1/2}}&  S(m,n,c)f(c) \\ & \approx\sum_{d|q}\mu(d)\sum_{\chi \text{ mod } q}\sum_{f\in \mathcal{B}(dq_1,\chi_1)}\frac{4\pi\sqrt{m'n'q_1^2}}{\cosh(\pi t_f)}\chi(m')\overline{\rho_f(m')}\rho_f(n'q_1^2)F(\chi)\widehat{H}_\infty(\pi_f) 
\end{align*}
under the constraint $X\geq \sqrt{mn}$. Here the inner sum ranges over the Maass forms of given level and nebentybus, $F$ is the Fourier transform of $f$, $\chi_1$ is the primitive character of $(\mathbb{Z}/q_1\mathbb{Z})^\times$ underlying $\chi$, and $\rho_f(m')$ with $m'=m/(m,q^\infty)$ denotes the Fourier coefficient normalised as in (3.2) loc. cit.\\ \noindent  Since their results rely on twisted multiplicativity of Kloosterman sums $S_\chi(m,n,c)$ it is not clear that their findings would generalise upon changing the neighborhood of $c$ (with respect to the $|\cdot|_p$-metric), nor that it would allow exponential sums other than the Kloosterman sums.
It is also worth pointing out that while their spectral expansion involves summing over automorphic forms with different nebentybus, the one appearing in Theorem \ref{thm:backward-formula} only runs over forms with trivial central character. \\
As a simple application we are able to obtain square-root cancellation (assuming Ramanujan-Petersson at all places) in sums of incomplete Kloosterman sums:

\begin{corollary}\label{cor:sum-of-kloosterman}
    Fix $m,n \in \mathbb{Z}$ such that $mn>0$ and $p$ a prime satisfying $(mn,p)=1$. For $q$ a sufficiently large $p$-power and $f_\infty \in \mathcal{S}(\mathbb{R}_+)$ and $X\geq \sqrt{mn}$, we have
    \begin{equation}
        \sum_{(c,p^\infty)=q}\frac{f_\infty\left(c/X\right)}{c^{1/2}}S(\overline{q}m,\overline{q}n,c_0)\ll_{f_\infty,\epsilon}X^{1/2+2\theta}(mnq)^\epsilon
    \end{equation}
    where $c=c_0 q$, the integer $\overline{q}$ is such that $\overline{q}q\equiv 1$ modulo $c_0$ and $\theta$ is any improvement towards the generalised Ramanujan-Petersson conjecture for arbitrary places.
\end{corollary}

\noindent Note that Corollary \ref{cor:sum-of-kloosterman} can be used to study exponential sums of Dedekind sums $S(d,c)$ which appear in the transformation laws of the Dedekind $\eta$-function (see \cite[(3.1)]{vardi-1987} for the transformation identity, also \cite{Rouse2015} for the connection between modular forms and $\eta$-quotients). This in turn can be used to obtain equidistribution results for particular sequences of Dedekind sums. See \cite{vardi-1987} for the equidistribution result on  $\langle rS(d,c)\rangle$ where $r \in \mathbb{R}_+$ with $c>0,0<d<c$ and $(d,c)=1$. Using identities relating Dedekind and Kloosterman sums we obtain
\begin{corollary}\label{cor:dedekind} Let $q$ be a $p$-prime power. Then
$$\sum_{c\leq X}\sum_{\substack{0<d<c_0 \\ (c_0,d)=1}}e(12\overline{q}S(d,c_0))\ll X^{1+2\theta}$$
In particular, the sequence
$\{\langle 12\overline{q}S(d,c_0) \rangle: d \text{ mod }  c_0, \ c\leq X, (c,p^\infty)=q  \}$
becomes equidistributed in $[0,1)$ as $X\to \infty$. 
\end{corollary}

\begin{remark}[Bessel functions in the literature]
   Applications of Bessel functions attached to representations of $\text{GL}_2$ to number theoretic questions come from the joint work of Baruch and Mao (\cite{BaruchMao2003},\cite{BaruchMao2007}), in relation to the Waldspurger correspondence and central values of automorphic L-functions. At the non-archimedean level, it is worth emphasazing that the Bessel function $i_{\pi,\psi}$ that they introduce is different from our Bessel function $j_{\pi,\psi}$, although $i_{\pi,\psi}$ and $j_{\pi,\psi}$ can be related to one another (see  \cite[Prop 4.3]{BaruchMao2003} and the paragraph immediately above). \\
   We also note that at the archimedean place representation theoretic Bessel functions have been used in \cite{Poincare-spectral} to re-prove the geometric-spectral formula \eqref{eq:backward-kuznetsov}. 

\end{remark} 

\begin{remark}[Some generalisations] Representation theoretic Bessel functions exist for any generic irreducible representation $\pi$ of $\text{GL}_2(\mathbb{Q}_p)$ regardless of central character. Results such as Theorem \ref{thm:forward-formula} and Theorem \ref{thm:backward-formula} should still apply for suitable choices of test functions $h_p$ and $f_p$. The opposite sign case $mn<0$ can also be addressed by taking instead the Bessel kernels
$$J_\pi^-(x)=\begin{cases}
        0 \hspace{3.3cm} \pi\simeq \mathcal{D}(\kappa) \\
        \pi i\frac{I_{2it}(x)-I_{-2it}(x)}{\sinh(\pi t)} \hspace{1.1cm} \pi \simeq \pi(\text{sgn}^\epsilon |\cdot|^{it},\text{sgn}^\epsilon |\cdot|^{-it})
    \end{cases}$$
Most of the local theory in this paper should extend to arbitrary number fields as well.
\end{remark}

\subsection{Notation and conventions}
\subsubsection{Groups, rings and matrices}\label{section:groups} Let $\mathbb{Q}$ denote the rational numbers and $\mathbb{Q}_v$ its completion with respect to the (archimedean or non-archimedean) place $v$. Let $\mathbb{A}$ and $\mathbb{A}_\text{fin}$ denote the ring of adeles and finite adeles of $\mathbb{Q}$ respectively, viewing the latter as a subring via diagonal embedding. Given a unital commutative ring $R$, we let $G(R):=\text{PGL}_2(R)$, and denote the unipotent and diagonal subgroups by
$$N(R)=\left\{ n(x):=\begin{pmatrix}
    1 & x \\ & 1
\end{pmatrix}: x \in R\right\} \qquad A(R)=\left\{ a(y):=\begin{pmatrix}
    y & \\ & 1
\end{pmatrix}: y \in R^\times\right\}$$
Let $B= NA$ the Borel subgroup of $G$. By Bruhat decomposition $G=B \sqcup NwB$ where $w=\begin{pmatrix}
    & -1 \\ 1 &
\end{pmatrix}$. When convenient, given a place $v$, we write $G_v=G(\mathbb{Q}_v)$. Let $K_\infty=\text{SO}_2(\mathbb{R})$ and $K_p=G(\mathbb{Z}_p)$, the standard maximal compact subgroups of $G(\mathbb{Q}_p)$ respectively. We further define
$$K_{0,p}(N)=\left \{ \begin{pmatrix}
    a & b \\ c& d
\end{pmatrix}\in G(\mathbb{Q}_p): c\equiv 0 \mod N \right\}$$
and globally we take  $K_0(N)=\prod_{p}K_{0,p}(N)$, a compact open subgroup of $G(\mathbb{A}_\text{fin})$. If we work at a finite place and the context is clear, we abuse notation and write $K_0(N)$ for $K_{0,p}(N)$. Finally we let $[G]=G(\mathbb{Q})\backslash G(\mathbb{A})$ and $[N]=N(\mathbb{Q})\backslash N(\mathbb{A})$.

\subsubsection{Measure normalisations} We give $\mathbb{R}$ the Lebesgue measure $dx$ and the measure $d^\times x=dx/|x|$ on $\mathbb{R}^\times$. On $\mathbb{Q}_p$ we choose the Haar measure $dx$ that gives $\mathbb{Z}_p$ measure $1$. We then define the measure $d^\times x$ on $\mathbb{Q}_p^\times$ by $d^\times x= \zeta_p(1)dx/|x|$, thus giving $\mathbb{Z}_p^\times$ volume $1$.\\
Let $dk$ be the probability Haar measure on $K_p$, and take measures on $A(\mathbb{Q}_v)$ and $N(\mathbb{Q}_v)$ induced by $d^\times x$ and $dx$ via isomorphism. These determine a Haar measure $dg$ on ${G}_{v}$ via the Iwasawa decomposition. \\
In the paper it will be useful to compare $dg$ against the measure $dg'$ pushed forward from $dn\times dxd^\times y$  on $N\times B\to G, (n,b)\mapsto nwb$. As discussed in \cite[\S 3.1.6]{Michel_2010}, we have the relation
$$d'g = \iota_v dg, \qquad \iota_v = \frac{\zeta_v(1)}{\zeta_v(2)}$$
where $\zeta_\infty(s)=\pi^{-s/2}\Gamma(s/2)$.

\subsubsection{Additive characters} We start by taking  $\psi = \prod_v \psi_v$ and $\psi_{\text{fin}}=\prod_{p<\infty}\psi_p$ where
$$\psi_v(x_v) = \begin{cases}
    e^{-2\pi i x_\infty} &  v=\infty ,\\
    e^{2\pi i\{x_p\}_p} &v=p
\end{cases}$$
where $\{x_p\}_p$ denotes the $p$-fractional part of $x_p$. All the non-trivial additive characters of $\mathbb{Q}\backslash \mathbb{A}$ have the form $\theta_\beta$ where
\begin{equation}\label{def:additive-char}
    \theta_\beta(x):=\psi(-\beta x)=\overline{\psi(\beta x)}
\end{equation}
for some $\beta \in \mathbb{Q^\times}$. When the context is clear, we identify $\theta_\beta$ with its restriction to a given place. \\
Given an additive local or global character $\psi$, we consider the space  $C_c^\infty(N\backslash G,\psi)$ consisting of smooth functions on $G$, compactly supported modulo $N$, satisfying the left transformation
\begin{equation}
    f(ng)=\psi(n)f(g) \qquad (n \in N, g \in G)
\end{equation}
We have an averaging map $A: C_c^\infty(G)\to C_c^\infty(N\backslash G,\psi)$ given by
\begin{equation}\label{def:averaging-map}
     f \mapsto \left(g \mapsto\int_{N}f(ng)\overline{\psi(n)}dn \right)
\end{equation}
See Lemma \ref{lemma:right-K(M)-invariance} for details.
\bigbreak

\subsection{Structure of the paper} We structure the paper as follows:
\begin{itemize}[itemsep=0em]
    \item Section 2 covers the necessary background material, including the basics of representation theory, the definition of Bessel functions, local-global relations for norms and Whittaker model, and the local Langlands correspondence for supercuspidal representations.
    \item Section 3 is dedicated to the integral transform $h_p \mapsto H_p$. We prove Theorem \ref{thm:forward-formula}, where issues of convergence are addressed by proving trivial bounds for $H_p$. Such bounds are then refined when $h_p$ is the indicator on a supercuspidal representation, to then derive Corollary \ref{cor:supercuspidal-selection}.
    A digression on the non-existence of newform projectors on the line $c(\pi)=1$ or $2$ and on the results mentioned in Remark \ref{remark:h-to-f} are also included.
    \item Section 4 provides an adelic framework of the key ideas in \cite{Poincare-spectral}, where we construct an adelic Poincaré series and prove its main properties.
    \item Section 5 is focused on the proof of Theorem \ref{thm:backward-formula} by an explicit computation of the spectral expansion of the aforementioned Poincaré series.
    \item In Section 6 we consider the transform $f\mapsto \widehat{H}$ in greater detail. We study the action of the Weyl element on a basis for the Kirillov model of generic $\pi \in \operatorname{PGL}_{2}(\mathbb{Q}_p)^\wedge$, providing explicit formulas for $c(\pi)\geq 2$  (and for some basis elements if $c(\pi )\leq 1$), obtaining bounds in the remaining cases. This is then applied to derive Corollary \ref{cor:sum-of-kloosterman} and Corollary \ref{cor:dedekind}.
\end{itemize}

\subsection*{Acknowledgements}
The author wishes to thank Ian Petrow for a careful reading of the draft of this paper and for providing many useful comments on its mathematical contents and exposition. The author is also grateful to Catinca Mujdei for a helpful discussion on the application of some results of this paper to shifted moments, and Chung-Hang Kwan for a discussion on deriving the geometric-spectral formula from the spectral-geometric one. This work was supported by the Engineering and Physical Sciences Research Council [grant numbers
EP/T517793/1, EP/W524335/1].

\subsection*{Statement on AI usage}
No generative AI was used at any point of the preparation of this paper.

\section{Background theory}
\subsection{Representation theory}\label{section:basics-on-local-theory}
\subsubsection{Basics on local theory} At the archimedean place, we denote by $\mathcal{D}(\kappa)$ the weight $\kappa$ discrete series of $G(\mathbb{R})$ as described in \cite[\S 11.7]{Traces}. Let $\chi=\text{sgn}^\epsilon |\cdot|^{it}$ a quasi-character of $\mathbb{R}^\times$, we construct the principal series $\pi(\chi,\chi^{-1})$ consisting of smooth functions $\phi:G(\mathbb{R})\to \mathbb{C}$ such that 
\begin{equation}\label{eq:principal-series-def}
    \phi(\begin{pmatrix}
    a & b \\ & d
\end{pmatrix}g)=\left|\frac{a}{d}\right|^{1/2}\chi(a)\chi^{-1}(d)\phi(g)
\end{equation}
We say $\pi=\pi(\chi,\chi^{-1})$ has parity $\varepsilon$ and spectral parameter $t$. We remark that we consider the archimedean representations to be $(\mathfrak{g},K)$-modules, but for sake of convenience we abuse terminology by referring to them as representations. We equip $G(\mathbb{R})^\wedge$ with the Plancherel measure $\widehat{\mu}$ given by
$$\widehat{\mu}(\pi)=\begin{cases}
    \frac{\kappa-1}{4\pi} & \pi\simeq \mathcal{D}(\kappa),\\
    \frac{1}{4\pi} t \tanh(\pi t)dt & \pi\simeq \pi(\text{sgn}^\varepsilon|\cdot|^{it},\text{sgn}^\varepsilon|\cdot|^{-it}), t \in i\mathbb{R},
\end{cases}$$
and is the zero measure otherwise. \bigbreak 
\noindent At finite places, for $\chi$ a quasi-character of $\mathbb{Q}_p^\times$, we construct the principal series $\pi(\chi,\chi^{-1})$ of $G(\mathbb{Q}_p)$ using \eqref{eq:principal-series-def}. If $\chi^2=|\cdot|^{\pm 1}$, we sometimes abuse notation and write $\pi(\chi,\chi^{-1})$ for the Steinberg representation $\chi St$. The unitary dual is equipped with a measure supported on tempered representations such that for $\pi$ square integrable $\widehat{\mu}(\pi)$ is equal to its formal degree, and if $c(\pi)=0$ then
$$d\widehat{\mu}(\pi)=\frac{(1+p^{-1})}{2\pi}\frac{(1-e^{2i\theta})(1-e^{-2i\theta})}{(1-p^{-1}e^{-2i\theta})(1-p^{-1}e^{2i\theta})}d\theta$$

\bigbreak

\noindent Let $\psi$ be a non-trivial additive character of $\mathbb{Q}_v$. Recall that an irreducible unitary generic representation of $G(\mathbb{Q}_v)$ admits a Whittaker model $\mathcal{W}(\pi,\psi)$ and a Kirillov model $\mathcal{K}(\pi,\psi)$, and that the map
$$\mathcal{W}\to \mathcal{K}, \quad W\mapsto (x\mapsto W(a(x)))$$
defines an isomorphism. \bigbreak

\noindent For $\pi$ a smooth irreducible unitary generic representation of $G(\mathbb{Q}_p)$, set 
\begin{equation}\label{def:euler-factors}
    \mathcal{L}_\pi(1)=\begin{cases}
        \frac{(1-p^{-2})}{(1-e^{2i\theta}p^{-1})(1-p^{-1})(1-e^{-2i\theta}p^{-1})} \text{ if } \pi\simeq \pi(|\cdot|^{i\theta/\log p},|\cdot|^{-i\theta/\log p}) , \\
        (1+p^{-1})^{-1} \hspace{2.65cm} \text{if } c(\pi)=1, \\
        (1-p^{-1}) \hspace{3.0cm} \text{if } c(\pi)\geq 2.
    \end{cases}
\end{equation}
where $\theta \in [0,\pi]$ or $\theta=i\tau\log p$ or $\pi+i\tau\log p$ with $\tau \in (0,1/2)$. Given a global standard (in the sense of \cite[\S 2.2.1]{Michel_2010}) automorphic representation $\pi\simeq \bigotimes_v'\pi_v$ with associated $L$-functions $\sum |\lambda_\pi(n)|^2 n^{-s}$ the leading Laurent coefficient $\mathcal{L}_\pi^*(1) $ at $s=1$ admits the factorisation
$$\mathcal{L}_\pi^*(1)=\prod_p\mathcal{L}_{\pi_p}(1)$$
where the product is regularised as in \cite[\S 4.1.5]{Michel_2010}.

\subsubsection{Bessel functions}\label{section:Bessel-functions} Working locally, fix a Whittaker functional $\varphi \mapsto W_{\varphi,\psi}(1)$ for some choice of non-trivial additive character $\psi$. For $x \in\mathbb{Q}_v^\times$, one may also define another choice of Whittaker functional given by
$$\varphi \mapsto \int_{\mathbb{Q}_v}W_\varphi(\begin{psmallmatrix}
    & x \\ -1 &
\end{psmallmatrix}\begin{psmallmatrix}
    1 & y \\ & 1
\end{psmallmatrix})\overline{\psi(y)}dy$$
We define the representation theoretic \emph{Bessel function} $j_{\pi,\psi}:\mathbb{Q}_v^\times \to \mathbb{C}$ to be the function obeying the relation
\begin{equation}\label{def:p-adic-bessel-function}
    \int_{\mathbb{Q}_v}W_\varphi(\begin{psmallmatrix}
    & x \\ -1 &
\end{psmallmatrix}\begin{psmallmatrix}
    1 & y \\ & 1
\end{psmallmatrix})\overline{\psi(y)}dy=j_{\pi,\psi}(x)W_\varphi(1)
\end{equation}
for all $x \in \mathbb{Q}_v^\times$. Its existence is a simple consequence of uniqueness of Whittaker functionals. At the archimedean place, for all $W\in \mathcal{W}(\pi,\psi)$, we have
\begin{equation}\label{eq:Bessel-fourier-inversion}
    W(\begin{psmallmatrix}
        & x \\ -1 &
\end{psmallmatrix})=\frac{1}{\zeta_\infty(1)}\int_{\mathbb{R}^\times}j_{\pi,\psi}(xy)W(\begin{psmallmatrix}
        y & \\ & 1
    \end{psmallmatrix})d^\times y
\end{equation}
(see eg \cite[\S 4]{Poincare-spectral}). At finite places instead if $W(a(x))\in C_c^\infty(\mathbb{Q}_p^\times)$, we have similarly
\begin{equation}\label{eq:fourier-inversion-bessel}
    W(\begin{psmallmatrix}
        & x \\ -1 &
\end{psmallmatrix})=\frac{1}{\zeta_p(1)}\int_{\mathbb{Q}_p^\times}j_{\pi,\psi}(xy)W(a(y))d^\times y
\end{equation}
(see \cite[Lemma 4.2]{Bessel}, noting that $\zeta_p(1)$ is absent by a different choice of normalisation). If $v=\infty$ it is possible to express $j_{\pi,\psi}$ in terms of standard Bessel functions (see e.g. \cite[Proposition 6.1]{Poincare-spectral}). For the weight $\kappa$ discrete series $\pi\simeq \mathcal{D}(\kappa)$ we have the relation
\begin{equation}\label{eq:archim-bessel-relation1}
    j_{\pi,\theta_n}(x)=\begin{cases}
        i^\kappa 2\pi|n|\sqrt{x}J_{\kappa-1}(4\pi|n|\sqrt{x}) \quad {x>0} \\
        0 \hspace{3.8cm} {x<0}
    \end{cases}
\end{equation}
whereas if $\pi$ is a principal series representation with spectral parameter $t$ (indipendent of parity) we have
\begin{equation}\label{eq:archim-bessel-relation2}
    j_{\pi,\theta_n}(x)=\frac{i\pi|n|\sqrt{x}}{\sinh(\pi t)}\begin{cases}
        J_{2it}(4\pi|n|\sqrt{x})-J_{-2it}(4\pi|n|\sqrt{x}) \qquad x>0 \\
        I_{2it}(4\pi|n|\sqrt{x})-I_{-2it}(4\pi|n|\sqrt{x}) \ \qquad x<0
    \end{cases}
\end{equation}
where $J_{\alpha}$ and $I_{\alpha}$ denote the standard Bessel functions.
 \bigbreak

\noindent Explicit expressions for $j_{\pi,\psi}$ exists at finite places as well. For instance if $\pi\simeq \pi(\chi,\chi^{-1})$, following \cite[\S 1.11]{Godement2018} we have
\begin{equation}\label{def:induced-model-bessel}
    j_{\pi,\psi}(x)=|x|^{1/2}\chi(x)\sum_{n\in \mathbb{Z}}\int_{v(z)=n}\psi(xz+z^{-1})\chi^2(z)d^\times z
\end{equation}
where the summation is to be understood as $|n|\leq N$ as $N\to \infty$. Similarly suppose $\pi$ is supercuspidal constructed from an admissible pair $(E/\mathbb{Q}_p,\beta)$ , where $E/\mathbb{Q}_p$ is quadratic and $\beta$ a character on $E^\times $ (see Section \ref{section:LLC} for details); a result of Baruch and Snitz \cite[Theorem 1.1]{Baruch2011} gives
\begin{equation}\label{eq:baruch-snitz}
    j_{\pi,\psi}(x)=\gamma\int_{E^x}\beta(z)\psi(\text{tr}(z))d_xz
\end{equation}
where $E^x$ is the subset of $E$ consisting of elements of norm $x$ with measure $d_x z$ and $\gamma$ is the Jacquet-Langlands constant attached to $E/\mathbb{Q}_p$. The author is not aware at present of similar explicit formulas for $j_{\pi,\psi}$ when $\pi$ is a supercuspidal representation that cannot constructed this way, though this can only occur at $p=2$ for small conductor exponents (see eg \cite[Thm 6.8]{generalised}). 

\subsubsection{Eisenstein series}\label{section:eisenstein-series}
For a character $\chi$ of $\mathbb{Q}^\times \backslash \mathbb{A}^\times$ we construct the global induced representation $\pi_s(\chi):=\bigotimes_v'\pi(\chi_v|\cdot|^s,\chi^{-1}|\cdot|^{-s})$. There is an isomorphism of global Hilbert spaces $\pi_0(\chi) \longrightarrow \pi_s(\chi)$ via the map $\phi \mapsto \phi_s$, where for $g=n(x)a(y)k \in G(\mathbb{A})$ with $x\in \mathbb{A}, y\in \mathbb{A}^\times, k\in K_\infty\times K$ and $\phi_s$ is given  by $\phi_s(g)=|a/d|^s\phi(g)$. For $\phi \in \pi_0(\chi)$ and $g \in G(\mathbb{A})$ we define the Eisenstein series $E(\phi_s,g)$ by the infinite sum
$$E(\phi_s,g)=\sum_{\gamma \in B(\mathbb{Q})\backslash G(\mathbb{Q})}\phi_s(\gamma g)$$
which converges absolutely for $\text{Re}(s)>1$. The Eisenstein series admits a meromorphic continuation to $s \in \mathbb{C}$.
Unless $\chi$ is quadratic and $s=1/2$, which are called the \emph{singular parameters}, we can equip $E(\phi_s,\cdot)$ with the Eisenstein norm
\begin{equation}||E(\phi_s,\cdot)||^2_{\text{Eis}}=\int_{K_\infty\times K}|\phi(k)|^2dk
\end{equation}
There exists an intertwiner $\pi_s(\chi)\to C^\infty([G])$ via $\phi_s \mapsto E(\phi_s,\cdot)$ (see \cite[\S 2.2.1]{Michel_2010}).

\subsubsection{Spectral decomposition}
Let $L^2([G])$ denote the Hilbert space of measurable, square integrable functions on $[G]$. The group $G(\mathbb{A})$ acts via the right-regular action $R(g)\varphi(x)=\varphi(xg)$. Following \cite[Thm. 6.2]{KL-kuznetsov}, the representation $R$ admits a decomposition into irreducible representations
\begin{equation}
    R=\bigoplus_{\pi \ \text{cuspidal}}\pi\  \oplus \ \bigoplus_{\chi^2=1}\pi_\chi \ \oplus \sum_{\substack{\chi \ \text{finite} \\ \text{order}}}\int_{-\infty}^{\infty}\pi_{it}(\chi)
\end{equation}
where $\pi_{\chi}$ are one dimensional representations spanned by the function $g \mapsto \chi(\det(g))$ for any quadratic Hecke character. Therefore, for any $\varphi \in L^2([G])$, we have the following spectral expansion (in the $L^2$-sense),
\begin{equation}\label{def:spectral-expansion}
    \varphi = \sum_{\pi \ \text{discrete}}F_\pi(\varphi)+\frac{1}{4\pi}\sum_{\substack{\chi \ \text{finite}\\ \text{order}}}\int_{-\infty}^{\infty}F_{\pi_{it}(\chi)}(\varphi)dt
\end{equation}
Here, if $\pi$ is discrete (ie. cuspidal or of the form $\pi_\chi$) we denote $F_\pi(\varphi)\in V_\pi$ the unique vector for which
\begin{equation}
    (\varphi,\varphi')=(F_\pi(\varphi),\varphi') \quad  \text{for all } \varphi'\in V_\pi.
\end{equation}
If instead $\pi$ is a global induced representation, and the function $\varphi$ is of compact support, we may form the inner product $(\varphi, E(\phi,\cdot))$ for $\phi \in \pi$. There exists a linear map $\Phi_\pi: C_c([G])\to V_\pi$ (which extends to $L^2([G])$) such that
\begin{equation}
    (\varphi, E(\phi,\cdot))=(\Phi_\pi(\varphi), \phi) \quad \text{for all } \phi \in V_\pi.
\end{equation}
We then take $F_\pi(\varphi)=E(\Phi_\pi(\varphi),\cdot)$.

\subsection{Dihedral local Langlands correspondence}\label{section:LLC}
Consider a quadratic field extension $E/\mathbb{Q}_p$ with residue $f$ and discriminant of valuation $d$. Given a character $\xi$ on $E^\times$, we say the pair $(E/\mathbb{Q}_p,\xi)$ is \emph{regular} if the character does not factor through the norm map, and we say it is \emph{admissible} if furthermore $E/\mathbb{Q}_p$ is unramified whenever $\xi|_{U_E(1)}$ factors through the norm map. Following  \cite[\S 39]{Bushnell2006} or \cite[\S 1.2]{newform-table}, it is possible to construct a supercuspidal representation $w_\xi$ for the group $\text{GL}_2(\mathbb{Q}_p)$ from a regular pair $(E/\mathbb{Q}_p,\xi)$. We say a represention of $\text{GL}_2(\mathbb{Q}_p)$ is \emph{dihedral} if it can be constructed this way.
\\
We consider the set of equivalence classes
$$\widetilde{\mathbb{P}}_2(\mathbb{Q}_p)=\{(E/F,\xi):\xi \text{ regular }, \xi \mid_{F^\times}=\chi_{E/F}\}/\sim_F$$
where $\chi_{E/F}$ is the non-trivial quadratic character on $F^\times /\text{Nm}(E^\times)$ and $(E/\mathbb{Q}_p, \xi)\sim_F (E'/\mathbb{Q}_p,\xi')$ if there exists a $F$-isomorphism $j: E\to E'$ such that $\xi=\xi'\circ j$.
For $p>2$, define
$$ \mathbb{P}_2(\mathbb{Q}_p)=\{(E/F,\xi) \in \widetilde{\mathbb{P}}_2(\mathbb{Q}_p) \text{ admissible}\}$$
so that there exists a bijection
\begin{equation}\label{eq:commutative-diagram}
    \mathbb{P}_2(\mathbb{Q}_p) \longrightarrow \{\text{supercuspidal }\pi \in G(\mathbb{Q}_p)\}, \quad (E,\xi)\mapsto \pi_\xi
\end{equation}
Similarly for $p=2$ we define instead
$$ \mathbb{P}_2(\mathbb{Q}_2)_{\geq 9}= \{(E/F,\xi) \in \widetilde{\mathbb{P}}_2(\mathbb{Q}_2): fc(\xi)+d\geq 9\}$$
then following \cite[Theorem 6.7]{generalised}, the diagram \eqref{eq:commutative-diagram} still holds upon replacing $\mathbb{P}_2(\mathbb{Q}_2) $ with $\mathbb{P}_2(\mathbb{Q}_2)_{\geq 9}$ as well as all supercuspidal representations with only those for which $c(\pi)\geq 9$ \bigbreak

\noindent Under this correspondance, one may express the $\epsilon$-factor of a dihedral supercuspidal representation $\pi$ constructed from $(E/\mathbb{Q}_p,\xi)$ as
\begin{equation}\label{eq:supercuspidal-epsilon-factor}
    \epsilon(s,\pi\otimes \chi,\psi)=\gamma\epsilon(s,\chi_E\xi,\psi_E)
\end{equation}
where $\chi_E=\chi\circ \text{Nm}_{E/F}$ and $\psi_E=\psi\circ \text{Tr}_{E/F}$, and $\gamma$ is the Jacquet-Langlands constant as given in \cite[\S34.3]{Bushnell2006}. The result is obtained by combining Thm 29.4 (3) and Corollary 30.4 loc.cit.

\subsection{Local-global relations}
We begin with the local setting, and describe how to relate Whittaker models with different additive characters. Let $\pi$ be a local generic representation. Fix a choice of additive characters $\psi$ and $\psi'$, then there exists $\beta \in \mathbb{Q}^\times_v$ such that $\psi(x)=\psi'(\beta x)$. Taking $W_{\varphi,\psi'} \in \mathcal{W}(\pi,\psi')$, the function
$$g \mapsto W_{\varphi,\psi'}(a(\beta)g)=\widetilde W_{\varphi,\psi}(g)$$
satisfies $\widetilde{W}_{\varphi,\psi}(ng)=\psi(n)\widetilde{W}_{\varphi,\psi}(g)$. Since left and right regular actions commute, the collection of all $\widetilde{W}_{\varphi,\psi}$ is another Whittaker model for $\pi$ with respect to $\psi$. Hence there is a non-zero constant $c(\pi,\psi,\psi')$ such that
\begin{equation}\label{def:c(pi,psi,psi')}
   c(\pi,\psi,\psi') W_{\varphi,\psi}(g)=W_{\varphi,\psi'}(a(\beta)g)
\end{equation}
\noindent Following \cite[Thm 6.2 (2)]{Petrow2023}, given an irreducible, infinite-dimensional, generic representation of $G(\mathbb{R})$, we consider the function in $\mathcal{K}(\pi,\theta_m)$ given by
\begin{equation}\label{def:minimal-vector}
    W_{\infty,m}(a(y))=\begin{cases}
        \text{sgn}(my)^\epsilon \left(\frac{\cosh(\pi t_\pi)}{\pi}\right)^{1/2}2\sqrt{|my|}K_{it_\pi}(2\pi|my|) \\
        \hspace{1.1cm}\left(\frac{(4\pi my)^\kappa}{\Gamma(\kappa)} \right)^{1/2}e^{-2\pi my}\ \delta(my>0)
    \end{cases}
\end{equation}
depending on whether $\pi$ is a principal series of parity $\epsilon$ and spectral parameter $t$, or a discrete series $\mathcal{D}(\kappa)$ respectively. The function $W_{\infty,m}$ is the unique minimal weight vector in $\mathcal{W}(\pi,\theta_m)$ satisfying the normalisation condition $||W_{\infty,m}||^2_{L^2(\mathbb{R}^\times)}=1$ and such that $W_{\infty,m}(1)>0$. \bigbreak 

\noindent Working globally, let $\pi$ be a standard generic automorphic representation of $G(\mathbb{A})$ and let $\varphi \in \pi$ smooth. For a given choice of (global) additive character $\psi$, we may define the Whittaker coefficient by
\begin{equation}\label{def:global-whittaker-coefficient}
    W_{\varphi,\psi}(g)=\int_{[N]}\varphi(n(x)g)\overline{\psi(x)}dx
\end{equation}
We call the vector space spanned by $W_{\varphi,\psi}$ for $\varphi \in \pi$ by  the \emph{global Whittaker model} of $\pi$. Recall that if $\pi \simeq \bigotimes'_v\pi_v$, we have an isomorphism $\mathcal{W}(\pi,\psi)\simeq \bigotimes'_v \mathcal{W}(\pi_v,\psi_v)$ so that if $\varphi\in \pi$ is a pure tensor, then
$$W_\varphi(g)=\prod_vW_{\varphi,v}(g_v).$$

\subsubsection{Inner products}\label{section:inner-products}
For a cuspidal automorphic representation (resp. a global principal series) $\pi$ and smooth $\varphi \in \pi$ (resp $\phi \in \pi$), we define the norm  
\begin{equation}\label{def:innerproduct}
    (\varphi,\varphi)_{L^2(\mathbb[G])}=\int_{[G]}|\varphi(g)|^2dg \qquad \left(\text{resp. }(\phi,\phi)=\int_{K_\infty \times K}|\phi(k)|^2dk \right)
\end{equation}
Suppose that $\pi\simeq \bigotimes_v'\pi_v$ with global Whittaker model $\mathcal{W}(\pi,\theta_m)$, we choose local inner products $(\cdot,\cdot)_v$ and local Whittaker models $\mathcal{W}(\pi_v,\theta_m)$ that are compatible in the following sense:
\begin{itemize}
    \item There are factorisations $\mathcal{W}(\pi,\theta_m)=\bigotimes_v\mathcal{W}(\pi_v,\theta_m)$ and $(\cdot,\cdot)=\prod_v(\cdot,\cdot)_v$ for pure tensors,
    \item If $\varphi_{0,\infty} \in \pi_\infty$ is the pre-image of $W_{\infty,m}$ given in \eqref{def:minimal-vector} under the Whittaker intertwiner $\pi_\infty\ \to \mathcal{W}(\pi_\infty,\theta_m)$, then $(\varphi_{0,\infty},\varphi_{0,\infty})_\infty=1$.
\end{itemize}
 We have the following: suppose that $\varphi=\varphi_{0,\infty}\otimes\bigotimes_{v<\infty}\varphi_v \in \pi$, where $\pi$ is a standard automorphic representation, that is right-$K_0(N)$-invariant for some $N \in \mathbb{N}$ and unit normalised with respect to the $L^2$ or Eisenstein norm, then
 \begin{equation}\label{eq:normalised-fourier-coefficients}
    \prod_p W_{\varphi_p,\psi}(1)=\frac{a_{u_\varphi}(m)}{\sqrt{|m|}}\delta_{m\in \mathbb{Z}}
\end{equation}
where $u_\varphi$ is the underlying Maass/holomorphic/Eisenstein form on $\Gamma_0(N)$ in the sense of  with normalised coefficients $a_{u_\varphi}(m)$. See \cite[Sections 2.2.1 and 2.3]{generalised} for a more detailed discussion. \bigbreak

\noindent Alternatively, given the factorisation $\mathcal{W}(\pi,\theta_m)=\bigotimes_v\mathcal{W}(\pi_v,\theta_m)$, we may define local regularised norms on $\mathcal{W}(\pi_v,\theta_m)$ by
\begin{equation}
    \langle W,W \rangle_v:=\frac{\int_{\mathbb{Q}_v^\times}W(a(y))\overline{W(a(y))}d^\times y}{\zeta_v(1)L(\pi_v,\operatorname{Ad},1)/\zeta_v(2)}
\end{equation}
and a global \emph{canonical norm} given by
\begin{equation}\label{def:canonical-norm}
    \langle W,W\rangle=\prod_v\langle W_v,W_v\rangle_v
\end{equation}
Following \cite[\S 2.2.2, \S 4.1.7, \S 4.4.2]{Michel_2010}  we obtain this result concerning inner products:

\begin{lemma} \label{lemma:norm-equality}
    Suppose $\pi$ is a standard generic automorphic representation, and let $\varphi_1,\varphi_2 \in \pi$ smooth. For an additive character $\psi$, there exists a constant $c_\psi(\pi)$ such that
    $$c_\psi(\pi)(\varphi_1, \varphi_2)=\langle W_{\varphi_1,\psi} , W_{\varphi_2,\psi}\rangle$$
    For all places $v$, there are local constants $c_{\psi,v}(\pi)$ such that $$c_{\psi,v}(\pi)(\varphi_{1v} ,\varphi_{2v})=\langle W_{\varphi_{1v},\psi}, W_{\varphi_{2v},\psi} \rangle$$
\end{lemma}

\begin{remark}
    In fact we could work uniquely with the regularised norms locally and the canonical norm globally. We choose to work with the inner products $(\cdot,\cdot)_v$ and global counterpart purely for sake of notation so that fewer constants appear in the computations. For this same reason, our definition of canonical norm \eqref{def:canonical-norm} also differs by a constant from that in \cite{Michel_2010}.
\end{remark}

\section{Spectral-geometric direction}
We briefly sketch the idea behind the derivation of Theorem \ref{thm:forward-formula}. Given $h_p \in \mathcal{S}(G(\mathbb{Q}_p)^\wedge)$ we construct a compactly supported test function $f_p \in C_c^\infty(G(\mathbb{Q}_p))$ such that for any generic $\pi \in G(\mathbb{Q}_p)^\wedge$ the operator $\pi(f_p)$ acts on the newform line by the scalar $h(\pi)$ and annihilates the orthocomplement. Using the relative trace approach of \cite{KL-kuznetsov}, we integrate an automorphic kernel underlying the test function $f_\mathbb{A}=f_\infty\otimes f_p \otimes\bigotimes_{\ell \neq p}1_{K_\ell}$ (see \S 3.3 and \S 3.4 for $f_\infty)$ to obtain Petersson-Kuznetsov type formulas where the off diagonal terms of the geometric side are expressed in terms of a double unipotent averaging of $f_\mathbb{A}$ (a Kloosterman sum as in \eqref{def:local-kloosterman-sum}). Using properties of the Bessel functions, we then recover the integral transform $H_\infty$. Overall we have a commutative diagram
\[\begin{tikzcd}
	{h_p} &&& {} & {H_p} \\
	&& {f_p}
	\arrow["{\operatorname{Integral \ transform}}\  \eqref{def:non-archim-H_transform}",from=1-1, to=1-5]
	\arrow["{\operatorname{Prop} \ref{prop:h_p-to-f_p}}"', dashed, from=1-1, to=2-3]
	\arrow[from=1-4, to=1-5]
	\arrow["{\operatorname{Kloosterman \ sum} \eqref{def:local-kloosterman-sum}}"', dashed, from=2-3, to=1-5]
\end{tikzcd}\]
but the information on $f_p$ needed is very minimal (essentially just that is compactly supported).

\subsection{Non-archimedean spectral decomposition}\label{section:non-arch spectral decomp}
We recall that for a $H$ a unimodular linear $p$-adic algebraic group (so separable and locally compact), we denote $H^\wedge$ the unitary generic dual. Letting $\Sigma$ be the Borel $\sigma$-algebra of $H^\wedge$ with respect to the Fell topology and $\mu$ a Haar measure on $H$ there exists a unique $\sigma$-finite measure $\hat{\mu}$ on $(H^\wedge,\Sigma)$, called the Plancherel measure, such that
\begin{equation}\label{eq:plancherel-formula}
\int_H|f(g)|^2d\mu = \int_{H^\wedge}||\pi(f)||_{\operatorname{HS}}^2d\widehat{\mu}
\end{equation}
for any $f \in L^1(H)\cap L^2(H)$ \cite[Thm 18.8.2]{Dixmier1969}. Here $||\cdot||_{\operatorname{HS}}$ denotes the Hilbert-Schmidt operator norm. \bigbreak

\noindent For a function $f\in C_c^\infty(H)$, we define $f^*_g:=\overline{f(h^{-1}g)}$,  the Plancherel theorem (\cite[\S 18.8.1]{Dixmier1969}) then asserts that
\begin{equation}\label{eq:plancherel-inversion}
    \overline{f(g)}=\int_{H^\wedge}\text{Tr}\pi(f_g^*)d\hat{\mu}(\pi)
\end{equation}
with the equality holding pointwise. We have the following mild generalisation of \cite[Prop. 4.1]{generalised}:

\begin{lemma}\label{lemma:non-archimedean-spectral-decomp}
    Let $f \in C_c^\infty(G(\mathbb{Q}_p))$. Suppose that for all $\pi \in G(\mathbb{Q}_p)^\wedge$ the operator $\pi(f):V_\pi \to V_\pi$ is self adjoint onto a finite dimensional subspace then we have the spectral expansion
    \begin{equation}
        f(g)=\int_{\pi \in G(\mathbb{Q}_p)^\wedge}\sum_{v \in \mathcal{B}_f(\pi)}\hat{f}_v(\pi)\overline{\Phi_{\pi,v}(g)}d\widehat{\mu}(\pi)
    \end{equation}
    where $\mathcal{B}_f(\pi)$ is any orthonormal eigenbasis for $\text{Im}(\pi(f))$, the function $\Phi_{\pi,v}(g)=\langle \pi(g)v,v\rangle$ is a matrix coefficient, and  $\hat{f_v}(\pi)$ is the eigenvalue associated to $v \in \mathcal{B}_f(\pi)$.
\end{lemma}

\begin{proof}
    Since $f \in C_c^\infty(G(\mathbb{Q}_p))$ we may apply Plancherel formula \eqref{eq:plancherel-formula}, moreover
    $$\text{Tr}(\pi(f_g^*))=\sum_{v\in \mathcal{B}(\pi)}\langle\pi(f_g^*)v,v \rangle$$ for $\mathcal{B}(\pi)$ any orthonormal basis extending $\mathcal{B}_f(\pi)$ (by hypothesis this is really a finite sum). Since $\pi(f_g^*)=\pi(g)\pi(f)^*$ we have $$\langle \pi(f_g^*)v,v\rangle = \langle v,\pi(f)\pi(g^{-1})v\rangle$$
    Exchanging summation and integration
    $$\text{Tr}(\pi(f_g^*))=\sum_{v\in \mathcal{B}(\pi)}\int_{G}\overline{f(h^{-1}g)}\Phi_{\pi,v}(h)dh=\sum_{v \in \mathcal{B}(\pi)}\int_{G}\overline{f(h)}\Phi_{\pi,v}(gh^{-1})dh$$
    and we note that
    \begin{align*}
        \int_{G}\overline{f(h)}\Phi_{\pi,v}(gh^{-1})dh=\int_{G}\overline{f(h)}\langle \pi(h^{-1})v, & \pi(g^{-1}v)\rangle dh\\
        & =\int_{G}\overline{f(h)}\sum_{w\in \mathcal{B}'(\pi)}\langle \pi(h^{-1})v,w\rangle \langle w,\pi(g^{-1})v \rangle dh
    \end{align*}
    where $\mathcal{B}'(\pi)$ is any orthonormal basis extending $\mathcal{B}_f(\pi)$ and respecting the decomposition $V_\pi = \text{Im}(\pi)\oplus \ker (\pi)$, possible by self-adjointedness of $\pi(f)$. Therefore
    $$\int_{G}\overline{f(h)}\Phi_{\pi,v}(gh^{-1})dh = \sum_{w\in \mathcal{B}'(\pi)}\langle v,\pi(f)w\rangle\langle w,\pi(g^{-1})v\rangle=\overline{\hat{f}_v(\pi)} \langle \pi(g)v,v\rangle$$
    since $w\in \mathcal{B}(\pi)$ is an eigenvector by assumption. The result follows upon taking conjugates and using \eqref{eq:plancherel-inversion}.
\end{proof}

\begin{proposition}\label{prop:h_p-to-f_p}
    Let $h_p\in \mathcal{S}(G(\mathbb{Q}_p)^\wedge)$. There exists a unique $f \in C_c^\infty(G(\mathbb{Q}_p))$ such that for all $\pi \in G(\mathbb{Q}_p)^\wedge$ the operator $\pi(f)$ acts on the newvector line by the scalar $h_p(\pi)$ and annihilates the orthocomplement, and so in particular $\operatorname{tr}(\pi(f))=h_p(\pi)$.
\end{proposition}
\begin{proof}
    We begin with existence, splitting the proof into several cases:
    \begin{itemize}[itemsep=0em]
        \item If $\sigma$  is supercuspidal, with normalised newvector $v_0$, let $f_\sigma=\frac{1}{||\Phi_{\sigma}||^2}\overline{\Phi}_{\sigma}$ where $\Phi_{\sigma}(g)=\langle{\sigma(g)v_0,v_0}\rangle$. By property of supercuspidal representations, the functions $\Phi_\sigma$ and $f_\sigma$ are compactly supported, locally constant. Moreover $\pi(f_\sigma)$ is a projector on $\mathbb{C}v_0$ if $\pi\simeq \sigma$ and it is the zero operator otherwise.
        \item Let $\chi$ be a character of $\mathbb{Z}_p^\times$ with $\chi^2$ non-trivial, from \cite[\S 7.2]{generalised} there exists $f_\chi \in C_c^\infty(G(\mathbb{Q}_p))$ such that $\pi(f_\chi)$ is a newform projector if $\pi\simeq \pi(\mu,\mu^{-1})$ with $\mu|_{\mathbb{Z}_p^\times}=\chi$ and $\pi(f_\chi)=0$ otherwise.
        \item Given a polynomial $P$ in one variable, by Satake isomorphism there exists a bi-$K_p$-invariant, compactly supported test function $f_p$ such that $\pi(f_p)$ acts by the scalar $P(2\cos \theta)$ on the newvector line of $\pi(|\cdot|^{i\theta/\log p},|\cdot|^{-i\theta/\log p})$ and annihilates the orthocomplement.
    \end{itemize}
    Existence then follows by taking linear combinations of such test functions. Uniqueness follows from Lemma \ref{lemma:non-archimedean-spectral-decomp}.
\end{proof}

\noindent Now assume that $\xi \in C_c^\infty(N(\mathbb{Q}_p)\backslash G(\mathbb{Q}_p),\psi)$, by the Whittaker-Plancherel formula (see \cite[Eq. (8)]{Jana2026}) we have
\begin{equation}\label{eq:whittaker-plancherel}
    \xi(g)=\int_{\pi\in G(\mathbb{Q}_p)^\wedge}\sum_{W\in \mathcal{B}(\pi)}W(g)\int_{N(\mathbb{Q}_p)\backslash G(\mathbb{Q}_p)}\xi(h)\overline{W(h)}dh d\widehat{\mu}(\pi)
\end{equation}
which converges absolutely. Here $\mathcal{B}(\pi)$ is a $L^2$-normalised orthonormal basis for $\mathcal{W}(\pi,\psi)$.

\begin{proposition}\label{prop:kloosterman-to-transform}
    Let $h_p \in \mathcal{S}(G(\mathbb{Q}_p)^\wedge)$ and $f_p\in C_c^\infty(G(\mathbb{Q}_p))$ the function attached to $h_p$ via Proposition \ref{prop:h_p-to-f_p}. For $(mn,p)=1$, the local Kloosterman sum of $f_p$ given in \eqref{def:local-kloosterman-sum} satisfies
    $$H_p(m,n,\mu)=H_p\left(\frac{mn}{\mu}\right).$$
\end{proposition}

\begin{proof}
    Using the bi-$A(\mathbb{Z}_p)$-invariance of $f \in C_c^\infty(G(\mathbb{Q}_p))$ we might rewrite 
$$H_p(m,n,\mu)=\int_{\mathbb{Q}_p}\int_{\mathbb{Q}_p}f(\begin{psmallmatrix}
        m  & -m t_1 \\ & 1 
    \end{psmallmatrix}\begin{psmallmatrix}
        & -\mu^{-1} \\ 1 & 
    \end{psmallmatrix}\begin{psmallmatrix}
        1 & nt_2 \\ & n
    \end{psmallmatrix})\psi_p(mt_1-nt_2)dt_1dt_2$$
    Letting $t_1'=mt_1, t_2'=nt_2$ and noting that
    $$\begin{pmatrix}
        m & -t_1 \\ & 1
    \end{pmatrix}\begin{pmatrix}
        & -\mu^{-1} \\ 1 & 
    \end{pmatrix}\begin{pmatrix}
        1 & t_2 \\ & n
    \end{pmatrix}=\begin{pmatrix}
        1 & -t_1 \\ & 1
    \end{pmatrix}\begin{pmatrix}
        & -mn \mu^{-1} \\ 1 & 
    \end{pmatrix}\begin{pmatrix}
        1 & t_2 \\ & 1
    \end{pmatrix}, $$
    we may rewrite
    $$H_p(m,n,\mu)=\iint_{\mathbb{Q}_p}f(\begin{psmallmatrix}
        1  & - t_1 \\ & 1 
    \end{psmallmatrix}\begin{psmallmatrix}
        & -mn\mu^{-1} \\ 1 & 
    \end{psmallmatrix}\begin{psmallmatrix}
        1 & t_2 \\ & 1
    \end{psmallmatrix})\psi_p(t_1-t_2)dt_1dt_2$$
    Denoting 
    $$ F(g)= \int_{\mathbb{Q}_p}f(\begin{psmallmatrix}
        1 & -t_1 \\ & 1
    \end{psmallmatrix}g)\psi_p(t_1)dt_1 \in C_c^\infty(N\backslash  G,\psi_p)$$
    we have by Whittaker-Plancherel expansion \eqref{eq:whittaker-plancherel} and unfolding the inner integral that
    \begin{align*}
        F(g)& =\int_{G(\mathbb{Q}_p)^\wedge}\sum_{v\in \mathcal{B}(\pi)}W_{v,\psi_p}(g)\int_{G(\mathbb{Q}_p)}f(g)\overline{W_{v,\psi_p}(g)}dg
         \\ & =\int_{G(\mathbb{Q}_p)^\wedge}h_p(\pi){W_{v_0,\psi_p}(g)}\overline{W_{v_0,\psi_p}(1)}d\widehat{\mu}(\pi)
    \end{align*}
    therefore 
    $$H_p(m,n,\mu)=\int_{\mathbb{Q}_p}\int_{G(\mathbb{Q}_p)^\wedge}h_p(\pi){W_{v_0,\psi_p}}(\begin{psmallmatrix}
        & -mn\mu^{-1}\\ 1 & 
    \end{psmallmatrix}\begin{psmallmatrix}
        1 & t \\ & 1
    \end{psmallmatrix})\psi_p(-t)\overline{W_{\varphi_0,\psi_p}(1)}d\widehat{\mu}(\pi) dt $$
    To exchange the order of integration, note that by the support of $f$ we may restrict the $t_2$-integral to some large compact set. Moreover the integral in equation \eqref{def:p-adic-bessel-function} stabilises upon restricting the integration to $|t_2|\leq p^N$ where $N$ depends only on the condutor of $\pi$, which is bounded by the restrictions on $h_p$ (see \cite[Lemma 4.1]{Bessel} and its proof). Therefore using \eqref{def:p-adic-bessel-function},
    $$H_p(m,n,\mu)=\int_{G(\mathbb{Q}_p)^\wedge}h_p(\pi)j_{\pi,\psi_p}(mn\mu^{-1})|W_{\varphi_0,\psi_p}(1)|^2d\widehat{\mu}(\pi)$$
    Noting that for $L^2(N\backslash G)$-normalised newforms (Lemma \ref{lemma:regularised-norm}), we have $$|W_{\varphi_0,\psi_p}(1)|^2=\zeta_p(1)^{-1}\mathcal{L}_\pi(1)^{-1}$$
    the theorem follows.
\end{proof}

\noindent We conclude this section with a result which may be of independent interest.
\begin{proposition}\label{prop:proj-nonexistence}
    Let $i \in \{1,2\}$. There exists no $f_{=i}\in C_c^\infty(G(\mathbb{Q}_p))$ such that $\pi(f_{=i})$ is the newform projector when $c(\pi)=i$ and is the zero operator otherwise.
\end{proposition}
\begin{proof}
    If such an $f_{=1}$ exists, then the spectral decomposition formula from Lemma \ref{lemma:non-archimedean-spectral-decomp} yields
    $$f_{=1}(g)=\frac{p-1}{2}\sum_{c(\pi)=1}\overline{\Phi_{\pi,v_0}(g)}$$
    where the sum runs over the Steinberg representation and its unramified twist. Letting $g=a(p^n)$ for $n$ even, and using the formulas for $\Phi_{\pi,v_0}$ given in \cite[Lemma 3.8]{Hu2018}, we obtain that $f_{=1}(a(p^n))\neq 0$ for infinitely many $n$, contradicting the assumption on the support.\\
    Similarly, let $h= \nu(p^2)1_{K_0(p^2)}-f_{=2}\in C_c^\infty(G(\mathbb{Q}_p))$ so that $\pi(h)$ is a projector onto $\pi^{K_0(p^2)}$ of $c(\pi)\leq 1$ and it is the zero operator otherwise. The spectral decomposition formula gives us 
    \begin{equation}\label{eq:matrix-coeff-1}
        h(g)=\int_{c(\pi)\leq 1}\sum_{v \in \mathcal{B}(\pi^{K_0(p^2)})}\overline{\Phi_{\pi,v}(g)}d\widehat{\mu}(\pi).
    \end{equation}
    By Atkin-Lenher theory, if $c(\pi)=0$ then $\pi^{K_0(p^2)}=\{v_1,v_p,v_{p^2}\}$ where $v_1$ is the newform for $\pi$ and $v_{p^j}=\pi(a(p^j))v_1$. Constructing an orthogonal basis $w_i=\sum_{j}A_{ij}v_{p^j}$ and using $G(\mathbb{Q}_p)$ invariance of the inner product we obtain
    \begin{equation}\label{eq:matrix-coeff-2}
        \Phi_{\pi,w_i}(a(p^n))=\sum_{j=0}^2 C_{ij}\Phi_{\pi,v_{1}}(a(p^{n-j}))
    \end{equation}
    for all large $n$ and some constants $C_{ij}$. Using the decomposition again gives
    $$1_{K_p}(a(p^n))=\int_{c(\pi)=0}\overline{\Phi_{\pi,v_1}(a(p^n))}d\widehat{\mu}(\pi)$$
    which is $0$ for $n\geq 1$. Now let $v_1$ be the newform in $St$, so that in the Kirillov model $v_1=|x|1_{\mathbb{Z}_p}(x)$, then denoting $\tilde{v}$ the orthocomplement of $v_1$ in $\pi^{K_0(p^2)}$, by Gram-Schmidt we obtain, up to scalar
    \begin{equation}\label{eq:matrix-coeff-3}
        \Phi_{\pi,\tilde{v}}(a(p^n))=(1+p^{-2})\Phi_{\pi,v_1}(a(p^n))-p^{-1}\Phi_{\pi,v_1}(a{p^{n-1}})-p^{-1}\Phi_{\pi,v_1}(a(p^{n+1}))=0
    \end{equation}
    using \cite[Lemma 3.8]{Hu2018}. The same result holds for the unramified twist of $St$. Combining equations \eqref{eq:matrix-coeff-1} to \eqref{eq:matrix-coeff-3} we reach the same contradiction as in $f_{=1}$ upon taking $n$ even.
\end{proof}

\subsection{Pre-trace formula}
The starting point for the Theorem \ref{thm:forward-formula} is an adelic pre-trace formula. Such formulas have originally appeared in work by Jacquet and Zagier, but for the purpose of this paper we make use of the Knigthly-Li's version \cite[equation (6.11)]{KL-kuznetsov}.

\begin{proposition}
    For any $f\in C_c^\infty(G(\mathbb{A}))$ and $(x,y) \in G(\mathbb{A})^2$ we have the decomposition
    $$K_{\operatorname{geom}}(x,y)=K_{\operatorname{cusp}}(x,y)+K_{\operatorname{cont}}(x,y)+K_{\operatorname{res}}(x,y)$$
    where 
    $$K_{\operatorname{geom}}(x,y)=\sum_{\gamma \in G(\mathbb{Q})}f(x^{-1}\gamma y)$$
    and
    \begin{gather*}
        K_{\operatorname{cusp}}(x,y)=\sum_{\pi \text{ cuspidal}}\sum_{\varphi \in\mathcal{B(\pi)}}R(f)\varphi(x)\overline{\varphi(y)}, \\
        K_{\operatorname{cusp}}(x,y)=\frac{1}{4\pi}\sum_{\chi \text{ finite order}}\sum_{\phi \in \mathcal{B}(\chi,\chi^{-1})}\int_{-\infty}^{\infty}E(\pi_{it}(f)\phi_{it},x)\overline{E(\phi_{it},y)}dt,\\
        K_{\operatorname{res}}(x,y)=\frac{3}{\pi}\sum_{\chi \text{ quadratic}}\chi(\det x)\overline{\chi(\det y)}\int_{G(\mathbb{A})}f(g)\chi(\det g)dg
    \end{gather*}
    Both sides converge absolutely and uniformly on compacta.
\end{proposition}

\begin{proof}
    See the discussion in \cite[\S 6]{KL-kuznetsov}.
\end{proof}

\subsection{Bruggeman/Kuznetsov formula derivation}\label{subsection-BK-formula} In this section we follow \cite{KL-kuznetsov}. Choose $f_\mathbb{A}=f_\infty\otimes f$ where $f_\infty$ is bi-$K_\infty$-invariant and supported on matrices of positive determinant and  $f=\bigotimes_{p<\infty}f_p$ is such that the function $f_p$ is chosen as in Proposition \ref{prop:h_p-to-f_p} and $f_q=1_{K_q}$ for $q\neq p$ . Following \cite[Prop. 3.6]{KL-kuznetsov}, there exists an isomorphism between $f_\infty$ as chosen above and $h$ in the even Paley-Weiner space. \bigbreak

\noindent We let $y_1,y_2>0$ and consider
$$I:=\frac{1}{\sqrt{y_1y_2}}\iint_{[N]^2}K_{\operatorname{geom}}(n_1\begin{pmatrix}
    y_1 & \\ & 1
\end{pmatrix},n_2\begin{pmatrix}
    y_2 & \\ & 1
\end{pmatrix})\overline{\theta_m(n_1)}\theta_n(n_2)dn_1dn_2$$
By compactness of $[N]$, we may apply Fubini's theorem to exchange the integral over $[N]^2$ with the sum defining $K_{\operatorname{cusp}}, K_{\operatorname{cts}}$ and $K_{\operatorname{res}}$, obtaining a decomposition
\begin{equation}\label{eq:kuznetsov-I-decomposition}
    I=I_{\operatorname{cusp}}+I_{\operatorname{cts}}+I_{\operatorname{res}}
\end{equation}
We further impose the constraint $w=my_1=ny_2$ and integrate both sides of \eqref{eq:kuznetsov-I-decomposition} in the interval $[0,\infty)$. On the spectral side 
$$\int_0^\infty I_{\operatorname{cusp}}(w)dw=\frac{1}{2}\sum_{\substack{\pi \ \text{cuspidal}\\ C(\pi)|p^\infty}}h(t_\pi)h_p(\pi_p)a_{u_\varphi}(m)\overline{a_{u_\varphi}(n)}$$
where we recall $\varphi$ denotes the $L^2$-normalised newform in $V_\pi$ (see proof of \cite[Prop. 7.5]{KL-kuznetsov}). Similarly
$$\int_{0}^\infty I_{\operatorname{cts}}(w)dw=\frac{1}{8\pi}\sum_{\chi: C(\chi)|p^{\infty}}h_p(\pi_{0}(\chi))\int_{-\infty}^{\infty}h(t)a_{u_{E(\pi_{it}(f)\phi_{it})}}(m)\overline{a_{u_{E(\phi_{it})}}(n)}dt$$
whereas the residual contribution vanishes. \bigbreak

\noindent To evaluate the geometric side, define for orbits $\delta \in N(\mathbb{Q})\backslash G(\mathbb{Q})/N(\mathbb{Q})$ the integral
\begin{equation}
    I_\delta(f_\mathbb{A})=\frac{1}{\sqrt{y_1y_2}}\int_{H_\delta(\mathbb{Q})\backslash H(\mathbb{A})}f_\mathbb{A}(\begin{psmallmatrix}
        y_1 & x_1 \\ & 1
    \end{psmallmatrix}^{-1}\delta \begin{psmallmatrix}
        y_2 & x_2 \\ & 1
    \end{psmallmatrix})\overline{\theta_m(x_1)}\theta_{n}(x_2)d(x_1,x_2)
\end{equation}
where $H=N\times N$, the group $H_\delta(\mathbb{Q})$ is the stabiliser of $\delta$ in $H(\mathbb{Q})$ acting on $G(\mathbb{Q})$ via $\gamma\cdot(x,y)=x^{-1}\gamma y$ and $d(x_1,x_2)$ is the quotient measure coming from $dn_1dn_2$.
Following \cite[\S 7.5]{KL-kuznetsov} we have
$$I=I_{\begin{psmallmatrix}
    n/m & \\ & 1
\end{psmallmatrix}}(f_\mathbb{A})+\sum_{\mu\in \mathbb{Q}^\times} I_{\begin{psmallmatrix}
    & -\mu \\  1 & 
\end{psmallmatrix}}(f_\mathbb{A}) $$
We call the first summand the \emph{first cell term}, and the remaining ones \emph{second cell terms}.

\subsubsection{First cell term}
We have to evaluate the integral
\begin{align}
    \int_0^\infty I_{\begin{psmallmatrix}
        n/m & \\ & 1
    \end{psmallmatrix}}& (f_\mathbb{A},w)dw \notag \\
    & =\delta_{m,n \in \mathbb{Z}}\delta_{mn>0}\frac{1}{2}\sqrt{\frac{n}{m}}\delta_\infty\int_{\mathbb{A}_{\operatorname{fin}}}f(\begin{pmatrix}
        n/m & t \\ & 1
    \end{pmatrix})\psi_{\operatorname{fin}}(-mt)dt
\end{align}
For $q\neq p$, since $f_q=1_{K_q}$ the local integral vanishes unless $n/m \in \mathbb{Z}_p$, in which case it evaluates to 1. By hypothesis $n/m \in \mathbb{Z}_p^\times$ as well therefore $m/n \in \{\pm1\}$. Using right-$A(\mathbb{Z}_p)$-invariance 

$$ \int_0^\infty I_{\begin{psmallmatrix}
        n/m & \\ & 1
    \end{psmallmatrix}} (f_\mathbb{A},w)dw =\frac{1}{2}\delta_{m=n}\delta_\infty\int_{\mathbb{A}_\text{fin}}f(\begin{pmatrix}
        1 & t \\ & 1
    \end{pmatrix})\psi_{\operatorname{fin}}(-mt)dt$$
so it remains to compute the non-archimedean integral.

\begin{lemma}\label{lemma:h-star}
    Suppose that $(nm,p)=1$. Let $h_p \in \mathcal{S}(G(\mathbb{Q}_p)^\wedge)$ and $f_p\in C_c^\infty(G(\mathbb{Q}_p))$ as in \ref{prop:h_p-to-f_p} be the associated function, then
    $$\int_{\mathbb{Q}_p}f_p(\begin{pmatrix}
        1 & t \\ & 1
    \end{pmatrix})\psi_{p}(mt)dt=\int_{G(\mathbb{Q}_p)^\wedge}\frac{h_p^*(\pi)}{\mathcal{L}_\pi(1)}d\widehat{\mu}(\pi)$$
    where $h_p^*(\pi)=h_p(\pi)$ if $c(\pi)\geq 1$ and for $c(\pi)=0$ is given by 
    $$h^*(\pi)=\frac{2}{\pi}\int_{-2}^{2}P(x)\sqrt{1-\frac{x^2}{4}}dx$$
\end{lemma}
\begin{proof}
    The case $c(\pi)\geq 1$ is addressed in the proof of \cite[Prop. 4.3]{generalised}. In the case $c(\pi)=0$, let $f_{p,d}$ the characteristic function on the subset of $M_{2\times 2}(\mathbb{Z}_p)$ with determinant in $p^d\mathbb{Z}_p^\times$, then by 
    \cite[Prop. 3.3]{KnigthlyLiPetersson},
    $$\int_{\mathbb{Q}_p}f_{p,d}(\begin{pmatrix}
        1 & t \\
        & 1
    \end{pmatrix})\psi_p(-mt)dt=\begin{cases}
        1 \quad \text{if } d=0\\
        0 \quad \text{otherwise}
    \end{cases}$$
 Moreover, for $\pi\simeq \pi(|\cdot|^{i\theta/\log p},|\cdot|^{-i\theta/\log p})$ we have $\pi(f_{p,d})=X_d(2\cos\theta)$ where $X_d$ is the Chebyshev polynomial of degree $d$ (see eg \cite[Prop. 29.8]{Traces}). Since Chebyshev polynomials are orthogonal with respect to the measure
 $\frac{2}{\pi}\sqrt{1-\frac{x^2}{4}}dx$ (see for instance \cite[Prop. 29.7]{Traces}) and any $f_p\in C_c^\infty(K_p\backslash G(\mathbb{Q}_p)/K_p)$ is a linear combination of $f_{p,d}$'s, the result follows upon applying \cite[Prop. 4.3]{generalised} again.
\end{proof}

\subsubsection{Second cell terms}
We may separate the archimedean and non-archimidean parts so that $I_\delta(f_\mathbb{A})=I_\delta(f_\infty)I_\delta(f)$ for $\delta=\begin{psmallmatrix}
    & -\mu \\ 1 & 
\end{psmallmatrix}$. Following the computations in \cite[\S 7.5.2]{KL-kuznetsov} we obtain
$$I_\delta(f_\infty)=\frac{1}{\sqrt{y_1y_2}}\iint_{\mathbb{R}^2}k(z_1,-\frac{\mu}{z_2})e(nx_2-mx_1)dx_1dx_2$$
where $z_j=x_j+iy_j$, the kernel is $k(z_1,z_2):=f(g_1^{-1}g_2)$ with $z_j=g_j\cdot i$, and
\begin{equation}
    I_\delta(f)=\iint_{\mathbb{A}_{\operatorname{fin}}^2}f(\begin{psmallmatrix}
        1 & -t_1 \\ & 1
    \end{psmallmatrix}\delta\begin{psmallmatrix}
        1 & t_2 \\ & 1
    \end{psmallmatrix})\psi_{\operatorname{fin}}(mt_1-nt_2)dt_1dt_2
\end{equation}
Relabeling $\mu \mapsto \mu^{-1}$ we see that $I_\delta(f)=H(m,n,x)$. We then require to exchange the order of summation and integration in
$$\int_0^\infty \sum_{\mu\in \mathbb{Q}_+}H(m,n,\mu)I_\delta(f_\infty,w)dw$$
This is justified by the following two results:
\begin{lemma}
    Suppose that $\{f_p\}_{p<\infty}$ is an admissible family of test functions with global Kloosterman sum $H(m,n,\mu)$. There exists a number $q\in \mathbb{Q}_+$ such that for all $\mu \in \mathbb{Q}_+\backslash q\mathbb{Z} $ and all $m,n \in \mathbb{Q}$ we have $H(m,n,\mu)=0$.
\end{lemma}
\begin{proof}
    Let $\mu=c^2x_0$ where $x_0\in \mathbb{Z}$ is square-free. Upon changing variables of integration
    $$H(m,n,\mu)=\iint_{\mathbb{A}_{\operatorname{fin}}^2}f(\begin{pmatrix}
        t_1 & \frac{-x_0^{-1}-t_1t_2}{c} \\ c & t_2
    \end{pmatrix})\psi_{\operatorname{fin}}\left(\frac{mt_1-nt_2}{c}\right)dt_1dt_2$$
    For $j$ squarefree, let $S_j\subset G(\mathbb{A}_{\operatorname{fin}})$ be a fundamental domain for the subset of $\text{supp} f$ whose matrices have determinant in $j\widehat{\mathbb{Z}}^\times$. Then $\text{supp} f$ is contained in as a finite union of such $S_j$, which are compact as they are the intersection of a compact and a closed set. The projection map $\pi:S_j \to \mathbb{A}_\text{fin}$ sending a matrix to its bottom left entry is continuous, hence $\pi(S_j)$ is also compact, which implies $H(m,n,\mu)$ vanishes unless $c$ has bounded denominator, which in turn implies $\mu$ has bounded denominator as $x_0$ can take only finitely many values.
\end{proof}

\begin{lemma}\label{lemma:kloosterman-properties}
    The generalised local Kloosterman sum $H_p(m,n,\mu)$ satisfies the following properties:
    \begin{enumerate}
        \item If $H_p(m,n,\mu)$ vanishes unless $m,n \in \mathbb{Z}_p$,
        \item If $H_p(m,n,\mu)$ is the Kloosterman sum for $f_p=1_{K_p}$ then given $\mu=x_0 p^{e_p(\mu)}$ we have
            $$H_p(m,n,\mu)=\begin{cases}
                S(\overline{x_0}m,n,p^{e_p(\mu)/2}) & e_p(\mu) \in 2\mathbb{N} \\
                0 & \text{otherwise}
            \end{cases}$$
        where $x_0\overline{x_0}\equiv 1$ modulo $p^{e_p(\mu)/2}$,
        \item There is a trivial bound of the form
        $$H_p(m,n,\mu)\ll ||f_p||_{L^\infty}|\mu|^{1/2}$$
        where the implicit constant depends on the support of $f_p$.
    \end{enumerate}
\end{lemma}

\begin{remark}
    In fact, more refined inequalities are available. If $h_p$ is the indicator on a principal series connected components this follows from \cite[Thm 3.8 (1)]{generalised}. We address in greater detail the supercuspidal case in Lemma \ref{lemma:supercuspida-trivial-bound}.
\end{remark}

\begin{proof} We address each point separately:
\begin{enumerate}
    \item  See the proof of \cite[Lemma 4.6 (3)]{generalised},
    \item Writing $\mu=x_0p^{e_p(x)}$ and change variables $t_1\mapsto x_0^{-1}t_1$ we obtain
    $$H_p(m,n,\mu)=\iint_{\mathbb{Q}_p^2}f_p\left(\begin{pmatrix}
        x_0^{-1}t_1 & x_0^{-1}(p^{-e_p(x)}-t_1t_2 \\
        1 & t_2
    \end{pmatrix}\right)\psi_p(mx_0^{-1}t_1-nt_2)dt_1dt_2$$
    We may write the matrix as 
    $$\begin{pmatrix}
        t_1 & p^{-e_p(x)}-t_1t_2 \\1 & t_2
    \end{pmatrix}\begin{pmatrix}
        x_0^{-1} & \\ & 1
    \end{pmatrix}$$
    which vanishes unless $e_p(x)$ is even. Using right-$A(\mathbb{Z}_p)$-invariance we have
    $$H_p(m,n,\mu)=S(\overline{x_0}m,n,p^{e_p(\mu)/2})$$
    where the last equality follows from \cite[Lemma 3.9 and eq. (3.15)]{generalised},
    \item By the Cartan decomposition we may assume that $\text{supp}f_p \subseteq K_p\begin{pmatrix}
        p^n & \\ & 1
    \end{pmatrix}K_p$,
    in this case there exists $\lambda \in \mathbb{Q}_p^\times$ for which the following system of equation holds
    \begin{equation}
        \begin{cases}
            2e_p(\lambda)-e_p(\mu)=n, \quad e_p(\lambda)+e_p(t_1)\geq 0, \\
            e_p(\lambda)+e_p(t_2)\geq 0, \qquad e_p(\lambda)\geq 0.
        \end{cases}
    \end{equation}
    so $e_p(\mu)$ and $n$ have the same parity. Set $\lambda =p^{\frac{n+e_p(\mu)}{2}}$, using center invariance and changing variables
    $$H_{p}(m,n,\mu)=\lambda^2\iint_{\mathbb{Q}_p^2}f\begin{pmatrix}
        t_1 & \frac{x_0^{-1}p^n-t_1t_2}{\lambda} \\ \lambda & t_2
    \end{pmatrix}\psi_p\left(\frac{mt_1-nt_2}{\lambda}\right)dt_1dt_2$$
    where we may restrict the integration to $(t_1,t_2)\in\mathbb{Z}_p$ such that $t_1t_2\equiv x_0^{-1}p^n \mod \lambda \mathbb{Z}_p$. Hence
    \begin{align*}
        H_{p}(m,n,\mu)& \leq \lambda^2||f||_{L^\infty}\int_{t_1 \in \mathbb{Z}_p }\int_{t_1t_2\equiv x_0^{-1}p^n \mod \lambda \mathbb{Z}_p} 1dt_1dt_2 \\
        & \ll \lambda^2||f||_{L^\infty}\frac{1}{\lambda}\ll ||f||_{L^\infty}p^{e_p(\mu)/2}
    \end{align*}
    which is the desired result.
\end{enumerate}
\end{proof}

\noindent Thus using the trivial bound in Lemma \ref{lemma:kloosterman-properties}, the Weyl bound for Kloosterman sums, and following the proof in \cite[Proposition 7.12]{KL-kuznetsov}, we obtain that
$$\sum_{\mu\in\mathbb{ Q}}H(m,n,\mu)\int_{0}^\infty I_{\delta}(f_\infty,w)dw$$
converges absolutely, and then the calculations in \cite[\S 7.5.2]{KL-kuznetsov} go through again with
$$\int_0^\infty I_\delta(f_\infty,w)dw=\frac{1}{4\sqrt{\mu}}\int_{-\infty}^{\infty} J_{2it}\left(4\pi\sqrt{mn\mu^{-1}}\right)\frac{h_\infty(t)t}{\cosh(\pi t)}dt=\frac{1}{2}\frac{H_\infty\left(4\pi \sqrt{mn\mu^{-1}}\right)}{\sqrt{\mu}}$$
where we abuse notation by denoting (in this section only)
$$H_\infty(x)=\frac{1}{4\pi}\int_{-\infty}^{\infty}h(t)t\tanh (\pi t)dt.$$
The extension to the larger class of test functions $h$ satisfying, for some $\delta>0$,
\begin{equation}\label{def:archim-test-function}
    \begin{cases}
        h(t) \text{ is even}, \\
        h(t) \text{ is holomorphic on the strip } |\text{Im}(t)|\leq 13,\\
        h(t)\ll (1+|t|)^{-2-\delta}
    \end{cases}
\end{equation}
is then carried out as in \cite[\S 8]{KL-kuznetsov}.

\subsection{Petersson formula derivation}\label{subsection:petersson} For convenience we address the case $\kappa\geq 4$ even only, though the $\kappa=2$ case can be addressed following \cite{DiScipio2026}. Let $v_0$ be an $L^2$-normalised lowest weight vector in the discrete series $\mathcal{D}(\kappa)$ and $\Phi_{\kappa,v_0}(g):=\langle\pi(g)v_0,v_0\rangle$. We replace the archimedean test function in the previous section \ref{subsection-BK-formula} with
$$f_{\infty}=\frac{1}{||\Phi_{\kappa, v_0}||^2_2}\overline{\Phi_{\kappa,v_0}}$$
or in explicit terms
$$f_\infty(g)=\begin{cases}
    \frac{\kappa-1}{4\pi}\frac{\det(g)^{\kappa/2}(2i)^\kappa}{(-b+c)+(a+d)i)^\kappa} \quad \text{ if }\det g>0 \\
    0 \hspace{3cm} \text{otherwise}
\end{cases}$$
Then the operator $\pi(f_\infty)$ is a projector on the newform line $v_0$ if $\pi_\infty\simeq \mathcal{D}(\kappa)$ and is identically zero otherwise. Integrating the kernel
$$I:=\iint_{[N]^2}K_{\operatorname{geom}}(n_1,n_2)\overline{\theta_{m}(n_1)}\theta_n(n_2)dn_1 dn_2$$
in two ways using the pre-trace formula we obtain $I=I_{\operatorname{cusp}}$. Using the same techniques of the previous section we obtain
\begin{equation}
    I_{\operatorname{cusp}}=\frac{(4\pi)^\kappa}{\Gamma(\kappa)}(mn)^{(\kappa-1)/2}e^{-2\pi(m+n)}\sum_{ \pi_\infty\simeq \mathcal{D}(\kappa)}h_p(\pi)a_{u_\varphi}(m)\overline{a_{u_\varphi(n)}}
\end{equation}
when $mn>0$ and $I_{\operatorname{cusp}}=0$ otherwise. Performing instead the geometric expansion into first and second cell terms we obtain
$$ I_{\begin{psmallmatrix}
    n/m & \\ & 1
\end{psmallmatrix}}=\delta_{m=n}\frac{(4\pi\sqrt{mn})^{\kappa-1}}{\Gamma(\kappa-1)}e^{-2\pi(m+n)}\int_{\mathbb{A}_\text{fin}}f(\begin{psmallmatrix}
    1 & t \\ & 1
\end{psmallmatrix})\psi_{\text{fin}}(-mt)dt$$
and
$$I_{\delta}(f_\mathbb{A})=\frac{(4\pi i)^{\kappa}(mn)^{(\kappa-1)/2}e^{-2\pi(m+n)}}{2\Gamma(\kappa-1)}\frac{H(m,n,x)}{\sqrt{x}}J_{\kappa-1}\left(4\pi\sqrt{\frac{mn}{x}}\right)$$
where $\delta=\begin{psmallmatrix}
     & -\mu^{-1} \\ 1 &
\end{psmallmatrix}$  and $H(m,n,\mu)$ is the global Kloosterman sum for the family $\{f_p\}_{p<\infty}$. Putting everything together we obtain
\begin{align}
    \sum_{\substack{\pi \text{ cuspidal} \\\pi_\infty\simeq \mathcal{D}(\kappa) \\C(\pi)|p^{\infty}}}h_p(\pi)a_{u_\varphi}(m)\overline{a_{u_\varphi}(n)}= &\delta_{m=n\in \mathbb{N}}\frac{\kappa-1}{4\pi}\int_{\mathbb{A}_{\text{fin}}}f(\begin{pmatrix}
        1 & t \\ & 1
    \end{pmatrix})\psi_{\text{fin}}(-mt)dt \notag \\
    &+ \frac{(\kappa-1)i^{-\kappa}}{2}\sum_{\mu\in \mathbb{Q}_+}\frac{H(m,n,\mu)}{\sqrt{\mu}}J_{\kappa-1}\left(4\pi\sqrt{\frac{mn}{\mu}}\right)
\end{align}

\subsubsection{Proof of Theorem \ref{thm:forward-formula}}
The theorem follows upon taking linear combinations of the Kuznetsov formula in Section \ref{subsection-BK-formula} and the Petersson formula in Section \ref{subsection:petersson} for finitely many choices of $\kappa$, and applying Proposition \ref{prop:kloosterman-to-transform}. The continuous term is given by
\begin{equation}\label{def:foward-thm-cts-term}
    \frac{1}{4\pi}\sum_{\chi: C(\chi)|p^\infty}h_p((\pi_0(\chi))_p)\int_{-\infty}^{\infty}h_\infty(t)a_{u_{E(\phi_{it})}}(m)\overline{a_{u_{E(\phi_{it})}}(n)}dt
\end{equation}
where
\begin{itemize}
    \item The sum ranges over finite order Hecke characters of $p$-power conductor,
    \item The representation $(\pi_0(\chi))_p$ is the $p$-th constituent of the global induced principal series as in Section \ref{section:eisenstein-series},
    \item $E(\phi_{it})$ is the Eisenstein series for the normalised newform $\phi \in V_{\pi_{it}(\chi)}$, as described in Section \ref{section:eisenstein-series}.
\end{itemize}

\begin{remark}\label{remark:recover-classical-formula}
    We briefly explain how to recover the classical Kuznetsov-Petersson formula for $\operatorname{SL}_2(\mathbb{Z})$ from Theorem \ref{thm:forward-formula}. Letting $h_p$ be the indicator on spherical representations, the spectral side is simply
    $$\sum_{\substack{\pi \ \text{cuspidal}\\ C(\pi)=1}}h_\infty(\pi_\infty)a_{u_\varphi}(m)\overline{a_{u_\varphi}(n)}+\frac{1}{4\pi}\int_{-\infty}^{\infty}h_\infty(t)a_{u_{E(\phi_{it})}}(m)\overline{a_{u_{E(\phi_{it})}}(n)}dt$$
    On the spectral side instead, the diagonal term evaluates to
    $$\delta_p=\frac{1}{2\pi}\int_{0}^{\pi}(1-e^{-2i\theta})(1-e^{2i\theta})d\theta=1.$$
    For the off-diagonal terms, note that from Proposition \ref{prop:kloosterman-to-transform} the integral transform $H_p(mn\mu^{-1})$ is supported on $e_p(\mu)\in 2\mathbb{N}$ and
    $$H_p\left(\frac{mn}{(cp^k)^{2}}\right)=S(\overline{c}m,\overline{c}n,p^k)$$
    appealing to $S(\overline{p^{2k}}m,n,c)=S(\overline{p^k}m,\overline{p^k}n,c)$ and twisted multiplicativity of Kloosterma sums we obtain
    $$\sum_{c\in \mathbb{N}}\frac{S(m,n,c)}{c}H_\infty\left(4\pi\frac{\sqrt{mn}}{c}\right)$$
\end{remark}

\subsection{Supercuspidal Kloosterman sums} 
In this section we prove some properties of the transform $H_p$ to deduce Corollary \ref{cor:supercuspidal-selection}.
We assume throughout that $\sigma$ is a dihedral supercuspidal representation of $\text{PGL}_2(\mathbb{Q}_p)$. If $p>2$, or $p=2$ and $c(\sigma)$ is odd, then all such representations are twist minimal, meaning $\min_{\chi}c(\sigma\otimes \chi)=c(\sigma)$. In $p=2$ and $c(\sigma)$ is even then \cite[Prop 6.37]{generalised} shows that $\min_\chi c(\sigma \otimes \chi)\geq c(\sigma)-2$.
Further let $h_p=1_{\pi\simeq \sigma}$, with transform $H_p=d_\sigma j_{\sigma,\psi}$ where $\psi$ is an additive character of trivial conductor.\\
Recall that we may associate a function $f\in C_c^\infty(G(\mathbb{Q}_p))$ to $h_p$ so that if $(mn,p)=1$ by Proposition \ref{prop:kloosterman-to-transform} we have the identity
$$H_p(m,n,\mu)=\iint_{\mathbb{Q}_p^2}f\begin{pmatrix}
    t_1 & \mu^{-1}-t_1t_2 \\ 1 & t_2 
\end{pmatrix}\psi(mt_1-nt_2)dt_2dt_2$$
If $(mn,p)\neq 1$ then we take the definition of $H_p(m,n,\mu)$ to be the integral above.

\begin{lemma}\label{lemma:supercuspida-trivial-bound}
    Fix $m,n \in \mathbb{Z}_p$ and $\mu \in \mathbb{Q}_p$. If $(mn,p)\neq 1$ then $H_p(m,n,\mu)=0$, whereas if $(mn,p)=1$ we have
    $$H_p(m,n,\mu)= \begin{cases}
        d_\sigma j_{\sigma, \psi_p}(mn\mu^{-1}) & e_p(\mu)\geq c(\sigma), \text{ and } \sigma \text{ twist min}, \\
        d_\sigma j_{\sigma, \psi_p}(mn\mu^{-1}) & e_p(\mu)\geq c(\sigma)-2, \text{ and }\sigma \text{ not twist min,} \\
        0 & \text{otherwise}.
    \end{cases}$$
    When non-vanishing, $H_p$ satisfies a trivial bound of the form
    $$ H_p(mnp^{-k})\leq d_\sigma p^{-\lfloor k/2\rfloor}\begin{cases}
        1 & \sigma \text{ is twist min}, \\
        p^2 & \sigma \text{ is not twist min}.
    \end{cases}$$
\end{lemma}

\begin{proof}
Suppose that $(mn,p)\neq 1$, we proceed as in the proof of Proposition \ref{prop:kloosterman-to-transform}, and obtain
    $$H_p(m,n,\mu) = d_\sigma\int_{t_2\in \mathbb{Q}_p}W_{\varphi_0,\theta_m}(\begin{psmallmatrix}
        & -\mu^{-1} \\ 1 &
    \end{psmallmatrix}\begin{psmallmatrix}
        1 & t_2 \\ & 1
    \end{psmallmatrix})\psi(-nt_2)dt_2 \overline{W_{\varphi_0,\theta_m}(1)}$$
    where $W_{\varphi_0}$ is the $L^2$-normalised newform for $\sigma$,
    and $W_{\varphi_0,\theta_m}(1)$ is up to a constant $W_{\varphi_0,\psi_p}(m)=0$ if $(p,m)\neq 1$. If $(p,m)=1$ instead, after changing Whittaker models from $\mathcal{W}(\sigma,\theta_m)$ to $\mathcal{W}(\sigma,\theta_n)$ as in \eqref{def:c(pi,psi,psi')}, to show vanishing of the integral above it suffices to evaluate 
    $$\int_{\mathbb{Q}_p}W_{\varphi_0,\theta_n}(\begin{psmallmatrix}
         & -m/(n\mu) \\ 1 
    \end{psmallmatrix}\begin{psmallmatrix}
        1 & t \\ & 1
    \end{psmallmatrix})\overline{\theta_m(t)}dt=j_{\sigma,\theta_n}(m/n\mu)W_{\varphi_0,\theta_n}(1)$$
    which again vanishes if $(p,n)\neq 1$.\\
    For the rest of the proof suppose that $(mn,p)=1$. Let $\chi$ be a character on $ \mathbb{Z}_p^\times$ and consider the Mellin transform
    $$\widehat{H}(\chi,k):=\int_{\mathbb{Z}_p^\times}H(mnup^{-k})\chi(u)d^\times u=d_\sigma \int_{\mathbb{Z}_p^\times}j_{\sigma,\psi_p}(mnup^{-k})\chi(u)d^\times u$$
    Combining the inversion formula \eqref{eq:Bessel-fourier-inversion} and \eqref{eq:supercuspidal-transform} we see that
    $$\widehat{H}(\chi,k)=d_\sigma \varepsilon(1/2,\sigma\otimes \chi^{-1},\psi)1_{\chi^{-1},-c(\sigma\otimes \chi^{-1})}(mnp^{-k})$$
    Thus $\widehat{H}(\chi,k)\neq 0$ if and only if $k=c(\sigma\times \chi^{-1})$. But for odd primes, all supercuspidal representations of trivial central character (in $\text{GL}_2(\mathbb{Q}_p)$) are twist minimal, meaning $c(\sigma\otimes \mu)\geq c(\sigma)$ for any quasi-character $\mu$ of $\mathbb{Q}_p^\times$, hence $k\geq c(\sigma)$. \\
    To obtain the upper bound, suppose that $\sigma$ is twist minimal, corresponding to $(E/\mathbb{Q}_p,\xi)$ under the LLC of Section \ref{section:LLC}. If $E/\mathbb{Q}_p$ is unramified then $c(\sigma\otimes \chi)=2c(\xi \chi_E)$ and so a necessary condition is $c(\chi_E)=c(\chi)\leq k/2$. Suppose instead $E/\mathbb{Q}_p$ is ramified, then $c(\pi\otimes \chi)=c(\xi\chi_E)+d$ where $d=1$ or $3$, the latter only if $p=2$. Since $c(\chi_E)=2c(\chi)-1$ we again require $c(\chi)\leq k/2$. The trivial bound then follows from Fourier inversion and the triangle inequality. The non-twist minimal case is similar.
\end{proof}

\begin{proof}[Proof of Corollary \ref{cor:non-vanishing-coeff}] Take  $h_\infty=1_{\mathcal{D}_\kappa}$ and $h_p=1_{\sigma}$, then Corollary \ref{cor:supercuspidal-selection} yields
$$\sum_{\substack{f \in S_\kappa^{\operatorname{new}}(N)\\ (\pi_f)_p\simeq\sigma}}|a_f(m)|^2=\frac{\kappa-1}{4\pi}\frac{d_\sigma}{\mathcal{L}_\sigma(1)}+\frac{\kappa-1}{2 i^\kappa} \sum_{\mu\in N\mathbb{N}}\frac{H^S(m,n,\mu)}{\sqrt{\mu}}H_p\left(\frac{mn}{\mu} \right)J_{\kappa-1}\left(\frac{4\pi m}{\sqrt{\mu}}\right)$$
Using the Weyl bound for $H^S(m,n,\mu)$ and the trivial bound in Lemma \ref{lemma:supercuspida-trivial-bound} for $H_p$, together with the fact that $J_\kappa(x)\ll_\kappa x^{-\kappa}$ for $x\leq 1$ yields
$$\sum_{\substack{f \in S_\kappa^{\operatorname{new}}(N)\\ (\pi_f)_p\simeq\sigma}}|a_f(m)|^2=\frac{\kappa-1}{4\pi}\frac{d_\sigma}{\mathcal{L}_\sigma(1)}+O_\kappa \left(d_\sigma \left(\frac{m^2}{N}\right)^{\kappa/2}\right)$$
    from which the result follows.
\end{proof}

\section{Geometric-spectral direction: preliminaries}
As in the previous section we provide a quick sketch on the derivation of Theorem \ref{thm:backward-formula}. Given a suitable test function $f\in C_c^\infty(N(\mathbb{A})\backslash G(\mathbb{A}),\psi)$, we construct an adelic Poincare series $P_f \in L^2([G])$, and evaluate its Whittaker functional in two ways: by unfolding the definition of $P_f$ and by performing a spectral expansion. This step is somewhat delicate due to convergence issues. Morally this is just an adelisation of the discussion in \cite{Poincare-spectral}, though the computations are more involved due to the local-global relations and the more flexible choice of $f$. 

\subsection{Adelic Poincaré series}
In this section we follow \cite[\S 2.3]{Poincare-spectral} to construct the adelic Poincaré series and prove the main properties which are later used to establish a version of the Kuznestov formula needed to derive Theorem \ref{thm:backward-formula}.
We begin by choosing a function $f \in C_c^\infty(N(\mathbb{A})\backslash G(\mathbb{A}),\psi)$, and defining the associated Poincaré series
\begin{equation}\label{def:Poincare-series}
    P_f(g)=\sum_{\gamma \in N(\mathbb{Q})\backslash G(\mathbb{Q})}f(\gamma g) \qquad (g \in G(\mathbb{A}))
\end{equation}

\begin{proposition}
    Let $f =\bigotimes_v f_v \in C_c^\infty(N(\mathbb{A})\backslash G(\mathbb{A}),\psi)$, then the Poincare series $P_f$ is compactly supported on $[G]$, and converges locally uniformly and absolutely on compacta.
\end{proposition}
\begin{proof}
    We begin by proving uniform and absolute convergence. Let $F$ be a fundamental domain for $\Gamma(1)\backslash \text{SL}_2(\mathbb{R})$ and $K=K_\infty\times K_\text{fin}$ so that $\mathcal{F}:=F\times K$ contains a fundamental domain for $G(\mathbb{Q})\backslash G(\mathbb{A})$ (see e.g. \cite[p.99]{KL-kuznetsov}). We prove absolute convergence on this set. Since $\mathbb{Q}$ has class number $1$, we may write $G(\mathbb{Q})=B(\mathbb{Q})\Gamma(1)$; for $g \in \mathcal{F}$ we have
    $$P_f(g)=\sum_{\eta \in \mathbb{Q}^\times}\sum_{\gamma \in \Gamma_\infty \backslash \Gamma(1)}f(a(\eta) \gamma g)$$
    At each finite place, since $\gamma g_p\in K_p$, we have that $|\eta|_p$ is bounded and in fact $\eta$ is a unit at almost all places. Therefore only finitely many $\eta$ thus give a non-zero contribution to $P_f$ and it suffices to show that
    $$P_f^{\eta}(g_\infty)=\sum_{\gamma \in \Gamma_\infty \backslash \Gamma(1)}f_\infty(a(\eta)\gamma g_\infty)$$
    converges absolutely.
    Taking $z=x+iy=g_\infty\cdot i$, for any $c>0$ the size of $|\{\gamma \in \Gamma_\infty\backslash \Gamma(1): \text{Im}(\gamma z)\geq c \}| $
    is bounded above independently of $z$ (see \cite[Lemma 2.10]{Spectralmethods}), therefore $P_f^\eta$ is locally a finite sum. In addition, for $y$ large, the height $H_\infty(a(\eta) \gamma g_\infty)$ is also large for all $\gamma \in \Gamma_\infty \backslash \Gamma (1)$, which by the compact support of $f_\infty$ implies that the sum is empty.
\end{proof}

\begin{lemma}\label{lemma:right-K(M)-invariance}
    Let $f \in C_c^\infty (N(\mathbb{A}_{\operatorname{fin}})\backslash G(\mathbb{A}_{\operatorname{fin}}),\psi)$, then there exists $M \in \mathbb{N}$ such that $f$ is right-$K(M)$-invariant.
\end{lemma}
\begin{proof}
    We prove that we can always express $f$ as a linear combination of averaging of indicator functions, in the sense of \eqref{def:averaging-map}. Let $g_i$ a choice of representatives for $N(\mathbb{A}_\text{fin})\backslash G(\mathbb{A}_\text{fin})$. Since $f$ is locally constant, we can find $M_i$ such that $f$ is constant on $g_iK(M_i)$. By compact support modulo $N$ we have
    \begin{equation}\label{eq:support-inclusion}
        \text{supp}(f)\subseteq \bigsqcup_{i=1}^kNg_iK(M_i).
    \end{equation}
    We may assume $f(g_i)\neq 0$ for $i=1,...,k$; but then the reverse inclusion in \eqref{eq:support-inclusion} also holds. Upon breaking the union into smaller double cosets if necessary
    $$\text{supp}(f) = \bigsqcup_{i=1}^{k'}Ng_iK(M)$$
    where $M=\max\{M_i,c(\psi)\}$. Therefore
    $$f=\sum_{i=1}^{k'}\frac{f(g_i)}{\tilde{1}_{g_i K(M)}(g_i)}\tilde{1}_{g_i K(M)},$$
    which is right-$K(M)$-invariant.
    Here $\tilde{1}_{g_iK(M)}$ denotes the image of $1_{g_iK(M)}$ under \eqref{def:averaging-map}. To show that the summation is well defined, it suffices to check that the denominator is non-vanishing: fixing $g=n(\xi)a(\eta)k$, we may use that $K(M)$ is normal in $K$ to rewrite this averaging as
    \begin{equation}
        \widetilde{1}_{gK(M)}(g)=\int_{N(\mathbb{A}_\text{fin})}1_{gK(M)g^{-1}}(n(x))\overline{\psi(x)}dx=\int_{N(\mathbb{A}_\text{fin})}1_{K(M)}(n(x\eta))\overline{\psi(x)}dx
    \end{equation}
    From the choice of $M$, the integral is non-vanishing.
\end{proof}

\begin{proposition}
    Fix $f \in C_c^\infty (N(\mathbb{A}\backslash G(\mathbb{A}),\psi)$. The Poincaré series $P_f$ is a smooth vector with respect to the $G(\mathbb{A})$ right-regular representation $R$.
\end{proposition}
\begin{proof}
    Let $\mathfrak{g}_\mathbb{C}$ denote the complexified Lie algebra of $G(\mathbb{R})$, and 
    take $X \in \mathcal{U}(\mathfrak{g}_\mathbb{C})$ (the universal enveloping algebra) and $g_\text{0}\in G(\mathbb{A}_\text{fin})$. By absolute convergence of the Poincaré series we have
    $$R(g_0)P_f = P_{R(g_0)f}\  \text{ and } \ R(X)P_f=P_{R(X)f}$$
    By Lemma \ref{lemma:right-K(M)-invariance}, the function $f$ has an open stabilizer in $G(\mathbb{A}_\text{fin})$ and therefore so does $P_f$. At $v=\infty$ instead, we note that $P_f$ is infinitely differentiable and such that $R(X)P_f \in L^2([G])$ for all $X$, since $R(X)f \in C_c^\infty(N(\mathbb{A})\backslash G(\mathbb{A}),\psi)$. Therefore smoothness with respect to $G(\mathbb{R})$ follows by \cite[Lemma 3.13]{Knapp}.
\end{proof}

\begin{lemma}\label{lemma:poincare-unfolding}
    Fix $f \in C_c^\infty(N(\mathbb{A})\backslash G(\mathbb{A}),\psi)$.
    If $\pi$ is automorphic cuspidal and and $\varphi \in V_\pi$ is a smooth vector, then
    \begin{equation}
        (P_f,\varphi)=\int_{N(\mathbb{A})\backslash G(\mathbb{A})}f(g)\overline{W_{\varphi,\psi}(g)}dg
    \end{equation}
    Similarly, if $\pi$ is a global induced representation, and  $\phi \in V_\pi$ is smooth, then
    \begin{equation}
        (P_f,E(\phi,\cdot))=\int_{N(\mathbb{A})\backslash G(\mathbb{A})}f(g)\overline{W_{E(\phi,\cdot),\psi}(g)}dg
    \end{equation}
\end{lemma}

\begin{proof}
    The lemma follows from a standard application of the unfolding principle, writing $P_f$ as an infinite sum
    \begin{equation}
        (P_f,\varphi)=\int_{[G]}P_f(g)\overline{\varphi(g)}dg=\int_{N(\mathbb{Q})\backslash G(\mathbb{A})}f(g)\overline{\varphi(g)}dg=\int_{N(\mathbb{A})\backslash G(\mathbb{A})}\int_{[N]}f(ng)\overline{\varphi(ng)}dndg
    \end{equation}
    and can conclude using $f(ng)=\psi(n)f(g)$. The derivation for global induced representations is similar.
\end{proof}

\subsection{The soft Kuznetsov formula}
For $\varphi \in  L^2([G])$ smooth and a global additive character $\psi$, we refer to the functional $\varphi \mapsto W_{\varphi,\psi}(1)$ induced from \eqref{def:global-whittaker-coefficient} as the \emph{standard Whittaker functional}. Let $f \in C_c^\infty(N(\mathbb{A})\backslash G(\mathbb{A}),\psi)$, choosing another character $\psi'$, we shall evaluate $W_{P_f,\psi'}$ in two ways, using the infinite summation definition of $P_f$ in \eqref{def:Poincare-series} to recover the geometric side, and the spectral expansion \eqref{def:spectral-expansion} to obtain the spectral side instead.
It is necessary to use a continuity argument to justify why the Whittaker functional commutes with the infinite summation in the spectral expansion.  \bigbreak

\noindent Suppose that $P_f$ is right-$K(M)$-invariant for some integer $M$ (using Lemma \ref{lemma:right-K(M)-invariance}). Consider the orbifold $\mathcal{M}:=G(\mathbb{Q})\backslash G(\mathbb{A})/K(M)$ and denote by $\mathcal{S}(\mathcal{M})$ the subspace of $L^2([G])$ consisting of smooth, right-$K(M)$-invariant vectors. We can equip $\mathcal{S}(\mathcal{M})$ with a finer Frechet topology induced by the seminorms $||R(X)\varphi||_{L^2(G(\mathbb{R}))}$ for $X \in \mathcal{U}(\mathfrak{g}_\mathbb{C})$. We state the following preliminary lemmas:

\begin{lemma}\label{lemma:whittaker-continuity}
The Whittaker functional $\varphi \mapsto W_{\varphi,\psi}(1)$ is continuous on the subspace of smooth vectors $\mathcal{S}(\mathcal{M})$ with respect to the finer topology.
\end{lemma}
\begin{proof}
    See the proof of  \cite[Theorem 1.1]{Poincare-spectral}. See also \cite[Corollary 2.5]{Wu_2014} and the remark immediately after.
\end{proof}

\begin{lemma}\label{lemma:dixmier}
    Let $\pi$ be a unitary representation of $G(\mathbb{R})$ with space of smooth vectors $V_\pi^\infty$. Given $v \in V_\pi^\infty$, there exist $f_1,...,f_n \in C_c^\infty(G(\mathbb{R}))$ and $v_1,...,v_n \in V_\pi$ such that
    $$v=\sum_{i=1}^n\pi(f_i)v_i$$
\end{lemma}
\begin{proof}
    This is a result due to Dixmier and Malliavin, see \cite{DixmierMalliavin1978}.
\end{proof}

\subsubsection{Geometric side} By the Bruhat decomposition $G=B \sqcup NwB$ we have
\begin{equation}
    P_f(g)=\sum_{\eta\in \mathbb{Q}^\times}f(a(\eta) g)+\sum_{\eta \in \mathbb{Q}^\times}\sum_{\xi \in \mathbb{Q}}f( w\ a(\eta) n(\xi)g)
\end{equation}
Let $\psi'=\theta_m$ and $\psi=\theta_n$. Then $W_{P_f,\theta_m}(g)$ consists of two terms, the \emph{diagonal} term and the \emph{off diagonal} term:
\begin{equation}
    \int_{[N]}P_f(xg)\overline{\theta_m(x)}dx=D+OD
\end{equation}
where
\begin{equation}\label{eq:geom-D}
    D =\sum_{\eta \in \mathbb{Q}^\times}\int_{[N]}f(a(\eta) n(x) g) \overline{\theta_m(x)}dx  
        =\sum_{\eta \in \mathbb{Q}^\times}f(a(\eta) g)\int_{[N]}\theta_n(\eta x)\overline{\theta_m(x)}dx= f(a(\beta) g)
\end{equation}
where $\beta=m/n$ using character orthogonality, and using $a(\eta)n(x)=n(\eta x)a(\eta)$. Moreover
\begin{equation}\label{eq:geom-OD}
    \begin{split}
        OD=\sum_{\eta\in \mathbb{Q^\times}}\sum_{\xi \in \mathbb{Q}}\int_{[N]}f( w\ a(\eta) n(x+\xi)g)\overline{\theta_m(x)}dx=\sum_{\eta \in \mathbb{Q}^\times}\int_{N(\mathbb{A})}f( wa(\eta)n(x)g)\overline{\theta_m(x)}dx.
    \end{split}
\end{equation}

\subsubsection{Spectral side} Fix a discrete automorphic representation or a global induced principal series $\pi$, as well as a test function $f\in C_c^\infty(N(\mathbb{A}\backslash G(\mathbb{A}),\theta_m)$.
Recall the projections $F_\pi(f):=F_\pi(P_f)$ as described in \S 2.1.4. We begin by proving the following results concerning projections of $P_f \in \mathcal{S}(\mathcal{M})$:
\begin{lemma}
    Suppose that $P_f \in \mathcal{S}(\mathcal{M})$. If $\pi$ is a discrete representation, the projection of $P_f$ onto $\pi$ satisfies $F_\pi(f) \in \mathcal{S}(\mathcal{M})$ as well. If  $\pi$ is a globally induced representation, then $F_\pi(f)$ is in $C^\infty(\mathcal{M})$.
\end{lemma}

\begin{proof}
    Let $k \in K(M)$ and any $\varphi \in V_\pi$. Since $\pi$ is unitary we have
    $$(P_f,\varphi)=(\pi(k)P_f,\varphi)=(\pi(k)F_\pi(f),\varphi)$$
    which implies $\pi(k)F_\pi(f)=F_\pi(f)$. Suppose $\pi$ is discrete, to prove smoothness with respect to $G(\mathbb{R})$, we appeal to Lemma \ref{lemma:dixmier}: we can find $f_1,...,f_n \in C_c^\infty(G(\mathbb{R}))$ and $v_1,...,v_n \in L^2([G])$ such that
    $$P_f= \sum_{i=1}^{n}R(f_i)v_i$$ By the linearity and intertwining properties of the spectral projections 
    $$F_\pi(f)=\sum_{i=1}^{n}\pi(f_i)F_\pi(v_i)$$
    and the result follows by applying Lemma \ref{lemma:dixmier} again with respect to the representation $\pi$. The case $\pi$ is not discrete is similar.
\end{proof}

\begin{proposition}[Soft Kuznetsov formula]\label{thm:soft-kuznietsov}
    Fix $m,n \in \mathbb{Q^\times}$ and $f \in C_c^\infty(N(\mathbb{A})\backslash G(\mathbb{A}),\theta_n)$. Expanding $W_{P_f,\theta_m}(1)$ in two ways, we obtain
    \begin{equation}
        \begin{split} 
        f(a(\beta))+ \sum_{\eta \in \mathbb{Q}^\times}\int_{N(\mathbb{A})}& f(wa(\eta)n(x))\overline{\theta_m(x)}dx = \\
        & =\sum_{\pi \in \Pi_c}W_{F_\pi(f),\theta_m}(1)+\frac{1}{4\pi}\sum_{\substack{\chi \ \text{finite}\\ \text{order}}}\int_{-\infty}^{\infty}W_{F_{\pi_{it}(\chi)},\theta_m}(1)dt
        \end{split}
    \end{equation}
    where $\beta=m/n$, the set $\Pi_c$ consists of cuspidal representations and the characters $\chi$ in the continuous term are finite order Hecke characters.
\end{proposition}

\begin{proof}
    Since $F_\pi(f)\in \mathcal{S}(\mathcal{M})$ (resp. $C^\infty(\mathcal{M})$) for all discrete representations (resp. global principal series) in the spectral expansion, it remains to show that the spectral decomposition converges with respect to $\mathcal{S}(\mathcal{M})$ to then appeal to the continuity of Lemma \ref{lemma:whittaker-continuity}. This is carried out in \cite[Proposition 1.4]{Poincare-spectral}. Note that the residue spectrum does not appear since it is annihilated upon applying the Whittaker functional and one dimensional representations are not generic. The proposition is then proved by putting together the spectral expansion \eqref{def:spectral-expansion} and \eqref{eq:geom-D},\eqref{eq:geom-OD}.
\end{proof}

\subsection{Archimedean inversion}
In this section we re-prove a generalised version of backwards Kuznetsov formula given in \eqref{eq:backward-kuznetsov} adapting the ideas in \cite{Poincare-spectral} to an adelic setting.
Starting with a function $f=\bigotimes_p 1_{K_p}$, and $f_\infty \in C_c^\infty(G(\mathbb{R}))$ with support in the open Bruhat cell, we construct a function $\tilde{f} \in C_{c}^\infty(N(\mathbb{A})\backslash G(\mathbb{A}),\theta_n)$ by averaging $f_\infty\otimes f$. \\ Explicitly, for any place finite place $v$, we define
\begin{equation}\label{def:left-char-test-function}
    \tilde{f}_{v}(g)=\int_{\mathbb{Q}_v}f_v(g^{-1}\begin{psmallmatrix}
        1 & t \\ & 1
    \end{psmallmatrix})\theta_n(t)dt
\end{equation}
this is the averaging of the function $g  \mapsto f_v(g^{-1})$ as described in \eqref{def:averaging-map}. The choice of having $g^{-1}$ in the argument allows to relate $\pi(f_v)$ with its adjoint operator, thus simplifying the computation of the spectral side.
At the archimedean place instead, we choose a test function supported on the open Bruhat cell: let $\eta \in C_c^\infty(\mathbb{R})$ and $\nu \in C_c^\infty(\mathbb{R}_+)$, one defines
\begin{equation}\label{def:inverted-test-function-arch}
    \widetilde{f}_\infty(g)=\begin{cases}
        \theta_n(x_1)\eta(x_2)\nu(y) & g = n(x_1)wn(x_2)a(y),\\
        \hspace{1.1cm} 0  & \text{otherwise.}
    \end{cases}
\end{equation}

\subsubsection{Recovering Kloosterman sums}
We compute the geometric side in the soft Kuznetsov formula of Proposition \eqref{thm:soft-kuznietsov} for our choice of $\widetilde{f}$ and show how Kloosterman sums arise in this setting.\\
Since $\widetilde{f}_\infty$ is not supported on $A(\mathbb{R})$, the diagonal term vanishes. By the support of $\widetilde{f}$ the summands in the off-diagonal term provide non-zero contribution only if they lie in $(\mathbb{Q}^\times)^2$. In addition $f_\infty$ is supported on matrices of positive determinant, we choose to reindex the sum over the $\eta$'s appearing in Proposition \ref{thm:soft-kuznietsov} as a sum over $c$ via $\eta=c^2$. Each summand splits into an archimedean and non-archimedean part.
The finite part is given by
\begin{equation}\label{eq:kloosterman-sums-appear}
    \int_{\mathbb{A}_\text{fin}^2}f(\begin{psmallmatrix}
        1 & -t_2 \\ & 1
    \end{psmallmatrix}\begin{psmallmatrix}
        & -c^{-2} \\ 1 & 
    \end{psmallmatrix}\begin{psmallmatrix}
        1 & t_1 \\ & 1
    \end{psmallmatrix})\psi_\text{fin}(-nt_1+mt_2)dt_2dt_1 = S(m,n,c)
\end{equation}
following \cite[Prop. 3.7]{KnigthlyLiPetersson}. The integral at infinity is instead given by
\begin{equation}\label{eq:int-at-infinity}
    \int_{\mathbb{R}}\widetilde{f}_\infty(wa(c^2)n(x))\overline{\theta_m(x)}dx
\end{equation}
Given that, by rearranging, 
$$\widetilde{f}_\infty(w\begin{pmatrix}
    y & \\ & 1
\end{pmatrix}\begin{pmatrix}
    1 & x \\ & 1
\end{pmatrix})=\widetilde{f}_\infty(w\begin{pmatrix}
    1 & yx \\ & 1
\end{pmatrix}\begin{pmatrix}
    y & \\ & 1
\end{pmatrix})=\eta(yx)\nu(y),$$
we may evaluate the integral \eqref{eq:int-at-infinity} as 
\begin{equation}\label{eq:backward-archim-appears}
    \frac{\nu(c^2)}{c^2}\int_{\mathbb{R}}\eta(x)\overline{\theta_m(x/c^2)}dx=\frac{\nu(c^2)}{c^2}\widehat{\eta}(\frac{1}{c^2})=\frac{A_{\eta,\nu}(c^{-2})}{c^2}, \qquad A_{\eta,\nu}(x):=\widehat{\eta}(x)\nu(x^{-1})
\end{equation}
where $\widehat{\eta}$ is the Fourier transform of $\eta$ with respect to ${\theta_m}$.

\bigbreak \noindent We now obtain an explicit expression for the spectral side of Theorem \ref{thm:soft-kuznietsov}.
We begin with the cuspidal contribution.
\begin{proposition}\label{prop:archimedean-transform}
    Fix $m,n\in \mathbb{Z}\backslash\{0\}$. Let $f_\mathbb{A}=f_\infty\otimes f \in C_c^\infty(G(\mathbb{A}))$, where $f_\infty$ is the test function supported on the open Bruhat cell as in \eqref{def:inverted-test-function-arch} and $f=\bigotimes_p1_{K_p}$.
    If $\pi$ is a generic standard automorphic representation, we have
    \begin{equation}
        W_{F_\pi(\widetilde{f}_\mathbb{A}),\theta_m}(1)=\frac{1}{\pi}\int_{\mathbb{R}^\times}A_{\eta,\nu}(y){j_{\pi,\theta_n}}(\beta y)d^\times y\sum_{\varphi \in \mathcal{B}_f(\pi)}a_{u_\varphi}(m)\overline{a_{u_{\varphi}}(n)}    \end{equation}
        where $\mathcal{B}_f(\pi)$ is an orthonormal basis consisting of pure tensors for $\text{Im}(\pi(f))$ that are minimal weight at $\infty$ and $\widetilde{f}_\mathbb{A}=A(f_\mathbb{A})$ is the periodisation given in \eqref{def:averaging-map}.
\end{proposition}

\begin{proof}
    Let $\varphi=\bigotimes_v \varphi_v$ be a pure tensor in $V_\pi$ such that each place is normalised by $(\varphi_v,\varphi_v)=1$, and take $W_\varphi$ its whittaker function. By Lemmas \ref{lemma:norm-equality} and \ref{lemma:poincare-unfolding} we have,
    \begin{equation}
        \langle W_{F_{\pi}(\widetilde{f}_\mathbb{A}),\theta_m},W_{\varphi,\theta_m}\rangle=c_{\theta_m}(\pi)\int_{N(\mathbb{A})\backslash G(\mathbb{A})}\widetilde{f}_\mathbb{A}(g)\overline{W_{\varphi,\theta_n}(g)}dg
    \end{equation}
    Using the definition \eqref{def:left-char-test-function} of $\widetilde f$ and the left action of unipotent elements on Whittaker models we obtain at the finite places,
    \begin{equation}\label{eq:unfolding-product}
        \int_{G(\mathbb{A}_{\operatorname{fin}})}f(g^{-1})\overline{W_{\varphi,\theta_n}(g)}dg=\overline{\int_{G(\mathbb{A}_{\operatorname{fin}})}f^*(g)W_{\varphi,\theta_n}(g)dg}=\prod_p \overline{\int_{G(\mathbb{Q}_p)}f_p^*(g)W_{\varphi_p,\theta_n}(g)dg}
    \end{equation}
    with $f^*(g):=\overline{f(g^{-1})}$. We now evaluate the integral on the right-hand side at each place. From Lemma \ref{lemma:norm-equality}, it follows that the Whittaker intertwiner  is bounded. By eg \cite[Proposition 10.10]{Traces}, we have  $$\pi(f^*_v)W_{\varphi,\theta_n}=W_{\pi(f_v^*)\varphi,\theta_n}.$$ 
    Suppose that $\pi(f_p)$ is a newform projector onto the newform line $\mathbb{C}\varphi_{0,p}$ of $\pi_p$ normalised so that $(\varphi_{0,p},\varphi_{0,p})=1$, then it is self adjoint operator and 
    $$\pi(f_p^*)\varphi_p=(\varphi_p,\varphi_{0,p})\varphi_{0,p}$$
    so that we have
    \begin{equation}\label{eq:proj-finite-place}
        W_{\pi(f_p^*)\varphi_p,\theta_m}(1)=(\varphi_p,\varphi_{0,p})W_{\varphi_{0,p},\theta_n}(1).
    \end{equation} 
    Working now in the archimedean setting, upon passing to the measure of the open Bruhat cell as in \S 1.2.2, we have
    \begin{equation}
        \int_{N(\mathbb{R})\backslash G(\mathbb{R})}f_{\infty}(g)\overline{W_{\varphi,\theta_n}(g)}dg = \iota_{\mathbb{R}}^{-1}\int_{\mathbb{R}}\int_{\mathbb{R}^\times}\eta(x)\nu(y)\overline{W_{\varphi,\theta_n}(w\begin{psmallmatrix}
            1 & x\\ & 1
        \end{psmallmatrix}\begin{psmallmatrix}
            y & \\ & 1
        \end{psmallmatrix})}dxd^\times y
    \end{equation}
    Given that $\beta=m/n$, the Whittaker function may be rewritten in terms of the Bessel function $j_{\pi,\theta_n}$ as follows
    \begin{equation}
        \begin{split}
            W_{\varphi,\theta_n}(w\begin{psmallmatrix}
                1 & x \\ & 1
            \end{psmallmatrix}\begin{psmallmatrix}
                y & \\ & 1
            \end{psmallmatrix})&= \int_{\mathbb{R}^\times}j_{\pi,\theta_n}(u)W_{\pi\begin{psmallmatrix}
                y & x \\ & 1
            \end{psmallmatrix}\varphi,\theta_n}\begin{psmallmatrix}
                u & \\ & 1
            \end{psmallmatrix}d^\times u \\
            & =\int_{\mathbb{R}^\times}\theta_n(ux)j_{\pi,\theta_n}(u)W_{\varphi,\theta_n}\begin{psmallmatrix}
                uy & \\ & 1
            \end{psmallmatrix}d^\times u \\
            & = \int_{\mathbb{R}^\times}\theta_m(\frac{ux}{ y})j_{\pi,\theta_n}(\frac{\beta u}{y})W_{\varphi,\theta_n}\begin{psmallmatrix}
                \beta u & \\ & 1
            \end{psmallmatrix}d^\times u
        \end{split}
    \end{equation}
    Recalling \eqref{def:c(pi,psi,psi')} and \eqref{def:minimal-vector}, it can be checked explicitly that $c(\pi_\infty,\theta_m,\theta_n)=1$, we may replace the Whittaker function on the right hand side with $W_{\varphi,\theta_n}(a(u))$. By the compact support of $\nu$, we may further exchange the order of integration to obtain
\begin{equation}\label{eq:proj-infinite-place}
    \begin{split}
    \int_{N(\mathbb{R})\backslash G(\mathbb{R})}f_{\infty}& (g)\overline{W_{\varphi,\theta_n}(g)}dg= \\&=\iota_\mathbb{R}^{-1}\int_{\mathbb{R}^\times}\left\{\int_{\mathbb{R}^\times}\hat{\eta}\left(\frac{u}{y}\right)\overline{j_{\pi,\theta_n}(\frac{\beta u}{y})}\nu(y)d^\times y \right\} \overline{W_{\varphi,\theta_m}(\begin{psmallmatrix}
        u & \\ & 1
    \end{psmallmatrix})}d^\times u
    \end{split}
\end{equation}
where $\widehat{\eta}$ is the Fourier transform of $\eta$ (see \eqref{def:inverted-test-function-arch}) with respect to $\theta_m$.
Suppose that $\pi(f)\neq 0$, using Lemma \ref{lemma:norm-equality}, and denoting by $\phi$ the $(\cdot,\cdot)$-normalised newform in $\pi$ we obtain from \eqref{eq:proj-finite-place} and \eqref{eq:proj-infinite-place} that
    \begin{align*}
        \langle W_{F_\pi(\widetilde{f}_\mathbb{A}),\theta_m},& W_{\varphi,\theta_m}\rangle=\\
        &=\frac{c_{\theta_m}(\pi)}{\iota_\mathbb{R}}\int_{\mathbb{R}^\times}\left\{\int_{\mathbb{R}^\times}\hat{\eta}\left(\frac{u}{y}\right)j_{\pi,\theta_n}(\frac{\beta u}{y})\nu(y)d^\times y \right\}\overline{W_{\varphi,\theta_m}(\begin{psmallmatrix}
        u & \\ & 1
    \end{psmallmatrix})}d^\times u(\phi,\varphi)_{\operatorname{fin}}\overline{W_{\phi,\theta_n}(1)}\\
        & =\frac{1}{\iota_\mathbb{R}}\left\langle \int_{\mathbb{R}^\times}A_{\eta,\nu}(uy){j_{\pi,\theta_n}}(\beta uy)d^\times y, W_{\varphi_\infty,\theta_m}\right\rangle\langle W_{\phi,\theta_m},W_{\varphi,\theta_m}\rangle_{\operatorname{fin}}\overline{W_{\phi,\theta_n}(1)}
    \end{align*}
    since the canonical norm relation implies $$(\phi,\varphi)_{\operatorname{fin}}=c_{\theta_m}(\pi)^{-1}(\zeta_\infty(1)L(1,\pi_\infty,\operatorname{Ad})/\zeta_\infty(2))^{-1}\langle W_{\phi},W_{\varphi}\rangle$$
    Therefore since the choice of $\varphi$ is arbitrary, and pure tensors span the vector space $V_\pi$,
    \begin{equation}
        W_{F_\pi(\widetilde{f}_\mathbb{A}),\theta_m}(1)=\frac{1}{\pi}\int_{\mathbb{R}^\times}A_{\eta,\nu}(y){j_{\pi,\theta_n}}(\beta y)d^\times y\ W_{\phi,\theta_m}(1)\overline{W_{\phi, \theta_n}(1)}
    \end{equation}
    the result then follows from \eqref{eq:normalised-fourier-coefficients}.
\end{proof}

\noindent \subsubsection{Archimedean transforms}\label{section:archimedean-transforms} Making use of \eqref{eq:normalised-fourier-coefficients} and changing variables $y\mapsto y^{-1}$ we have 
\begin{equation}\label{eq:W(1)-archimedean}
        W_{F_\pi(\widetilde{f}_\mathbb{A}),\theta_m}(1)= \frac{a_{u_\varphi}(m)\overline{a_{u_\varphi}(n)}}{\pi\sqrt{|mn|}}\times\int_{\mathbb{R^\times}}A_{\eta,\nu}(y){j_{\pi,\theta_n}}(\beta y)d^\times y
    \end{equation}
    if $\pi$ is an automorphic representation for which $\pi(f)\neq 0$, and $W_{F_\pi(\widetilde{f}_\mathbb{A}),\theta_m}(1)=0$ otherwise.
Recall the expressions \eqref{eq:archim-bessel-relation1}, \eqref{eq:archim-bessel-relation2}, and assume $mn>0$. If $\pi_\infty$ is a principal series of spectral parameter $t$, the integral transform in \eqref{eq:W(1)-archimedean} becomes
\begin{equation*}\begin{split}\frac{i\pi}{\sinh(\pi t)}& \int_{0}^{\infty}A_{\eta,\nu}(x)\sqrt{mnx} (J_{2it}(4\pi\sqrt{mnx}) -J_{-2it}(4\pi\sqrt{mnx})\frac{dx}{x} \\
& = \frac{2i\pi \sqrt{mn}}{\sinh(\pi t)}\int_{0}^{\infty}uA_{\eta,\nu}(u^2) (J_{2it}(4\pi u\sqrt{mn})-J_{-2it}(4\pi u \sqrt{mn}))\frac{du}{u}
\end{split}\end{equation*}
If instead $\pi_\infty\simeq \mathcal{D}(\kappa)$, the integral transform becomes
$$4\pi \sqrt{mn} \ i^\kappa \int_{0}^\infty u A_{\eta,\nu}(u^2) \ J_{\kappa-1}(4\pi\sqrt{mn}u)\frac{du}{u}$$
for weight $\kappa$ discrete series. The case $mn<0$ is similar, and we can see that the holomorphic terms vanish as expected. This suggests that really one should think of $g(x):=xA_{\eta,\nu}(x^2)$ as the archimedean test function appearing on the two sides of the trace formula, and since $\eta,\nu$ are arbitrary, any choice of $g$ in $C_c^\infty(\mathbb{R}_+)$ is admissible.

\subsubsection{Proof of equation \eqref{eq:backward-kuznetsov}}
On the geometric side, combine Proposition \ref{thm:soft-kuznietsov} with \eqref{eq:backward-archim-appears} and \eqref{eq:kloosterman-sums-appear} to obtain
\begin{equation}
    \sum_{c\in \mathbb{N}}\frac{A_{\eta,\nu}(c^{-2})}{c^2}S(m,n,c)
\end{equation}
Using the discussion in Section \ref{section:archimedean-transforms}, the contribution from holomorphic forms is given by
\begin{equation}
    \mathcal{H}=\sum_{\substack{\kappa>0\\\kappa \ \text{even}}}\sum_{\substack{\pi:\pi_\infty\simeq D(\kappa)\\C(\pi)=1}}a_{u_\varphi}(m)\overline{a_{u_\varphi}(n)}\frac{2}{\pi}\int_{0}^{\infty}\frac{A((u/4\pi\sqrt{mn})^2)u}{4\pi\sqrt{mn}}\cdot 2\pi i^{\kappa}J_{\kappa-1}(u)\frac{d^\times u}{u}
\end{equation}
whereas the contribution from Maass forms is given by
\begin{align}
    \mathcal{M}=\sum_{\substack{\pi \text{ Maass}\\ C(\pi)=1}}a_{u_\varphi}(m)\overline{a_{u_\varphi}(n)}\frac{2}{\pi}\int_0^\infty  \frac{A((u/4\pi\sqrt{mn})^2)u}{4\pi\sqrt{mn}}\cdot \frac{i\pi(J_{2it}(u)-J_{-2it}(u))}{\sinh(\pi t)}d^\times u
\end{align}
and the continuous contribution is similar:
\begin{equation}
    \mathcal{E}=\frac{1}{4\pi}\int_{-\infty}^{\infty}\frac{2}{\pi}(\dots)a_{u_{E(\phi_{it})}}(m)\overline{a_{u_{E(\phi_{it})}}(n)}dt
\end{equation}
where $(\dots)$ denotes the same integral transform as the Maass forms case. Taking $f_\infty(c)=A_{\nu,\eta}((c/4\pi\sqrt{mn})^2)(c/4\pi\sqrt{mn})$ then gives the result.

\section{Geometric-spectral direction: main results}
We now make slight modification to the arguments in the previous section in order to obtain a backwards formula at an arbitrary finite place $p$. On the geometric side, we weigh the generalised Kloosterman sums by some function $f_p$ in $C_c^\infty(\mathbb{Q}_p^\times)$, and on the spectral side we obtain a summation over all possible irreducible unitary generic representations at $G(\mathbb{Q}_p)$ weighted by an integral transform of $f_p$. We thus develop a pleasant analogy between the archimedean and the non-archimedean places, since in both cases we are summing over all possible representations of $G(\mathbb{Q}_v)$ which occur as local components of standard automorphic representations, and the transforms involved have a representation theoretic interpretation. 

\subsection{Derivation of the formula}
The heart of the derivation is analogous to the one given in the proof of Proposition \ref{prop:archimedean-transform}: we shall apply Bessel inversion to express $W_{\varphi,\psi}(w)$ in terms of $j_{\pi,\psi}$ and $W_{\varphi,\psi}(a(y))$. 

\noindent Define a complex valued function $\widetilde f_p$ on $G(\mathbb{Q}_p)$, which mirrors the choice of archimedean test function,
\begin{equation}\label{def:inverted-testfunction-nonarch}
    \widetilde f_p(g)=\begin{cases}
        \theta_n(x_1)\eta(x_2)\nu(y) & g = \begin{pmatrix}
            1  & x_1 \\ & 1
        \end{pmatrix}w \begin{pmatrix}
            1 & x_2 \\ & 1
        \end{pmatrix}\begin{pmatrix}
            y & \\ & 1
        \end{pmatrix} \\
        \hspace{1.2cm} 0 & \text{otherwise}
    \end{cases}
\end{equation}
where we choose Schwarz functions $\eta \in C_c^\infty(\mathbb{Q}_p)$ and $\nu \in C_c^\infty(\mathbb{Q}_p^\times)$. Again we have that $\widetilde f_p \in C_c^\infty(N(\mathbb{Q}_p)\backslash G(\mathbb{Q}_p),\theta_n)$ and we may define
$$A_{\eta,\nu}(\mu^{-1}):=|\mu|_p\int_{\mathbb{Q}_p}f_p(w a(\mu)n(x))\overline{\psi_m(x)}dx=\nu(\mu)\widehat{\eta}\left(\frac{1}{\mu}\right)$$
which shall appear on the geometric side. 

\bigbreak \noindent We construct $\widetilde f_\mathbb{A} =\bigotimes_v\widetilde{f_v}$ as follows: at the places $\infty$ and $p$, we take $\widetilde{f}_\infty$ and $\widetilde{f}_p$ as in \eqref{def:inverted-test-function-arch} and \eqref{def:inverted-testfunction-nonarch}. For all finite places $\ell\neq p$, we consider $\{f_\ell\}_{\ell\neq p}$ an admissible family of test functions, and take $\widetilde f_\ell$ to be the averaging in \eqref{def:left-char-test-function}.
We obtain the following result concerning the spectral side.
\begin{proposition}\label{prop:whittaker-functional-finite-inversion} Let $\widetilde{f}_\mathbb{A}$ the test function defined above. Fix $m,n \in \mathbb{Z}$ with $mn>0$. Let $\pi=\bigotimes_v'\pi_v$ a standard generic automorphic representation, then
    \begin{equation}
    \begin{split}
        W_{F_\pi(\widetilde{f}_\mathbb{A}),\theta_m}(1)=\frac{d(\pi_p,\theta_m,\theta_n)}{\pi \zeta_p(1)}\int_{\mathbb{R}^\times}A_{\eta,\nu}(y)&{j_{\pi,\theta_n}}(\beta y)d^\times y \int_{\mathbb{Q}_p^\times}A_{\eta,\nu}(y)  \overline{j_{\pi,\theta_n}(\beta y)}d^\times y\\
        &\times \prod_{\ell \neq p}\sum_{\varphi \in \mathcal{B}_{f_\ell}(\pi_\ell)}W_{\varphi,\theta_m}(1)\overline{W_{\varphi,\theta_n}(1)}
    \end{split}
    \end{equation}
    where for $\ell\neq p$, the basis $\mathcal{B}_{f_p}(\pi)$ consists of orthonormal, $(\cdot,\cdot)_p$-normalised vectors spanning $Im(\pi(f_p))$, and
    $$d(\pi_p,\theta_m,\theta_n)=c_{\theta_m,p}(\pi)\overline{c(\pi_p,\theta_m,\theta_n)}L(\pi_p,\operatorname{Ad},1)$$
\end{proposition}

\begin{proof} Again we choose an arbitrary $\varphi = \bigotimes_v \varphi_v$, this time with the additional requirement that $W_{\varphi_p,\theta_m}(a(y)) \in C_c^\infty(\mathbb{Q}_p^\times)$. The span of such $\varphi$ is dense in $V_\pi$.
Unfold as in \eqref{eq:unfolding-product}, by the definition of $\widetilde f_{p}$ we have
    \begin{equation}
        \int_{N(\mathbb{Q}_p)\backslash G(\mathbb{Q}_p)}f_{p}(g)\overline{W_{\varphi,\theta_n}(g)}dg = \iota_p^{-1}\int_{\mathbb{Q}_p}\int_{\mathbb{Q}_p^\times}\eta(x)\nu(y)\overline{W_{\varphi,\theta_n}(w\begin{psmallmatrix}
            1 & x\\ & 1
        \end{psmallmatrix}\begin{psmallmatrix}
            y & \\ & 1
        \end{psmallmatrix})}dxd^\times y
    \end{equation}
    Rewriting the Whittaker function in terms of Bessel functions via \eqref{eq:fourier-inversion-bessel} we obtain 
    \begin{equation}
            \zeta_p(1)W_{\varphi,\theta_n}(w\begin{psmallmatrix}
                1 & x \\ & 1
            \end{psmallmatrix}\begin{psmallmatrix}
                y & \\ & 1
            \end{psmallmatrix})=
             \int_{\mathbb{Q}^\times_p}\theta_m(\frac{ux}{ y})j_{\pi,\theta_n}(\frac{\beta u}{y})W_{\varphi,\theta_n}\begin{psmallmatrix}
                \beta u & \\ & 1
            \end{psmallmatrix}d^\times u
    \end{equation}
    \noindent By the compact support of $\eta$ and $\nu$, we may exchange the order of integration to obtain
\begin{equation}\label{eq:pre-canonical-non-archim}
    \begin{split}
    \int_{N(\mathbb{Q}_p)\backslash G(\mathbb{Q}_p)}f_{p}& (g)\overline{W_{\varphi,\theta_n}(g)}dg= \\&=\frac{\overline{c(\pi_p,\theta_m,\theta_n)}}{\iota_p\zeta_p(1)}\int_{\mathbb{Q}_p^\times}\left\{\int_{\mathbb{Q}_p^\times}\hat{\eta}\left(\frac{u}{y}\right)\overline{j_{\pi,\theta_n}}(\frac{\beta u}{y})\nu(y)d^\times y \right\} \overline{W_{\varphi,\theta_m}(\begin{psmallmatrix}
        u & \\ & 1
    \end{psmallmatrix})}d^\times u
    \end{split}
\end{equation}
Considering the rescaling by $L(1,\pi_p,\operatorname{Ad})\zeta_p(1)/\zeta_p(2)$ appearing in the regularised norm $\langle\cdot,\cdot\rangle_p$, equation \eqref{eq:pre-canonical-non-archim} becomes
$$\frac{L(1,\pi_p,\operatorname{Ad})\overline{c(\pi_p,\theta_m,\theta_n)}}{\zeta_p(1)}\left\langle \int_{\mathbb{Q}_p^\times}\widehat{\eta}(uy^{-1})\overline{j_{\pi_p,\theta_n}(\beta u y^{-1})}\nu(y)d^\times y,W_{\varphi_q,\theta_m} \right\rangle_p$$
The rest of the proof follows as before through an application of Lemma \ref{lemma:norm-equality} and computations at the other finite and archimedean place analogous to those in Proposition \ref{prop:archimedean-transform}.
\end{proof}

\subsection{Computation of the constants}
In this section we compute the constants $c_{\theta_m}(\pi_v)$ and $c(\pi_v,\theta_m,\theta_n)$ for $\pi_v$ an irreducible unitary generic representation at a finite place.
We shall now assume that $p\nmid mn$, in particular $\theta_m$ and $\theta_n$ have conductor $\mathbb{Z}_p$. \bigbreak
\noindent We start with computing $c_{\theta_m}(\pi_p)$. Let $\varphi_{0}$ denote a $( \cdot, \cdot)$-normalised newform for $\pi_p$. Then
\begin{equation}
    c_{\theta_m}(\pi_p)=\langle W_{\varphi_0,\theta_m}, W_{\varphi_0,\theta_m}\rangle ^2=|W_{\varphi_0,\theta_m}(1)|^2\left\langle\frac{W_{\varphi_0,\theta_m}}{W_{\varphi_0,\theta_m}(1)}, \frac{W_{\varphi_0,\theta_m}}{W_{\varphi_0,\theta_m}(1)}\right\rangle^2
\end{equation}

\begin{lemma}\label{lemma:regularised-norm} Let $\pi$ be an irreducible, generic representation of $G(\mathbb{Q}_p)$ with newform $W\in \mathcal{W}(\pi,\psi)$ (with $\psi$ having trivial conductor) normalised so that $W(1)=1$, then
$$\langle W,W\rangle=\frac{\zeta_p(2)\mathcal{L}_\pi(1)}{L(1,\pi,\operatorname{Ad})}$$
\end{lemma}
\begin{proof}
    It suffices to compute expressions for $||W||^2_{L^2(\mathbb{Q}_p^\times)}$ in terms of $\mathcal{L}_\pi(1)$. Suppose that $c(\pi)\geq 2$, then by \cite[p. 23]{newform-table} it is easy to see $||W||^2_{L^2(\mathbb{Q}_p^\times)}=1=\mathcal{L}_\pi(1)\zeta_p(1)$. Similarly if $c(\pi)=1$, we have $$||W||^2_{L^2(\mathbb{Q}_p^\times)}=\zeta_p(2)=\zeta_p(1)\mathcal{L}_\pi(1)$$
    Lastly, if $\pi$ is a spherical representation, following \cite[\S3.4.1]{Michel_2010} we have
    \begin{align*}
        ||W||^2_{L^2(\mathbb{Q}_p^\times)}& =\zeta_p(1)L(1,\pi,\operatorname{Ad})/\zeta_p(2)\\ &=\zeta_p(1)\frac{1-p^{-2}}{(1-e^{2i\theta}p^{-1})(1-p^{-1})(1-e^{-2i\theta}p^{-1})}=\zeta_p(1)\mathcal{L}_\pi(1)
    \end{align*}
    Rescaling by the regularising constant yields the result.
\end{proof}

\noindent Similarly we may write
\begin{equation*}
    W_{\varphi_0,\theta_n }(1)\frac{W_{\varphi_0,\theta_n}(\begin{psmallmatrix}
        \beta & \\ & 1
    \end{psmallmatrix})}{W_{\varphi_0,\theta_n}(1)}=c(\pi,\theta_m,\theta_n)W_{\varphi_0,\theta_m}(1)
\end{equation*}
the fraction on the left hand side depends only on $v(\beta)$ and so takes value 1. This implies
\begin{equation}
    c(\pi,\theta_m,\theta_n)=\frac{W_{\varphi_0,\theta_n}(1)}{W_{\varphi_0,\theta_m}(1)}
\end{equation}
thus combining Lemma \ref{lemma:regularised-norm} we obtain
\begin{equation}\label{eq:d(pi_q,theta_m,theta_n)}
    d(\pi_p,\theta_m,\theta_n)=\zeta_p(2)\mathcal{L}_{\pi_q}(1)W_{\varphi_0,\theta_m}(1)\overline{W_{\varphi_0{,\theta_n}}(1)}
\end{equation}
\bigbreak

\noindent Using the expression for $d(\pi_p,\theta_m,\theta_n)$ above and Proposition \ref{prop:whittaker-functional-finite-inversion}, given some standard generic automorphic $\pi$ we obtain
\begin{align}\label{eq:recovering-coefficient}
    \frac{d(\pi_p,\theta_m,\theta_n)}{\zeta_p(1)}\prod_{\ell \neq p}\sum_{\varphi \in \mathcal{B}_{f_\ell}(\pi_\ell)}W_{\varphi,\theta_m}(1) & \overline{W_{\varphi,\theta_n}(1)}\notag \\&=\frac{\mathcal{L}_{\pi_p}(1)\zeta_p(2)}{\sqrt{|mn|}\zeta_p(1)}\sum_{\varphi \in \mathcal{B}_f(\pi)}a_{u_{\varphi}}(m)\overline{a_{u_{\varphi}}(n)}
\end{align}

\subsection{Proof of theorem \ref{thm:backward-formula}}
The proof now follows from simple book-keeping, combining Propositions \ref{thm:soft-kuznietsov} and \ref{prop:whittaker-functional-finite-inversion}, the proof of Proposition \ref{prop:archimedean-transform} as well as equation \eqref{eq:recovering-coefficient} and the discussion in Section \ref{section:archimedean-transforms}.
In particular, we have a continuous term of the form
\begin{equation}\label{eq:backward-formula-cts-term}
    \mathcal{E}=\frac{1}{4\pi}\sum_{\chi}\sum_{\phi \in \mathcal{B}(\pi_0(\chi),\mathcal{F})}\int_{-\infty}^{\infty}\widehat{H}_{\infty}(\pi_\infty)\widehat{H}_p({\pi_p})a_{u_{E(\phi_{it})}}(m)\overline{a_{u_{E(\phi_{it})}}(n)}dt
\end{equation}
where the outer sum ranges over finite order Hecke characters.

\section{Non-archimedean transforms}

\subsection{The Bessel transform}
Let $\psi$ be a non-trivial additive character on $\mathbb{Q}_p$ with conductor $\mathbb{Z}_p$. Throughout this section only, we take $dx$ to be the Haar measure that is self dual to $\psi$, giving $\mathbb{Z}_p$ volume 1, and we also take $d^\times x=dx/|x|$. Given a character $\chi$ of $\mathbb{Z}_p^\times$ and $n \in \mathbb{Z}$ we define
\begin{equation}
    1_{\chi,n}(x)=
    \begin{cases}
        \chi(u) & \text{if} \quad x=up^n \ (u\in \mathbb{Z}_p^\times, n\in \mathbb{Z}) \\
         \quad  0 & \text{otherwise}
    \end{cases}
\end{equation}
which form a spanning set for $C_c^\infty(\mathbb{Q}_p^\times)$. Fixing some $\beta \in \mathbb{Q}_p^\times$, for $\pi \in {G(\mathbb{Q}_p)}^\wedge$ we define the \emph{Bessel transform} of $1_{\chi,n}$ at $\pi$ by
\begin{equation}
    B_\pi(1_{\chi,n})=\int_{\mathbb{Q}_p^\times}1_{\chi,n}(u)j_{\pi,\psi}(\beta u)d^\times u
\end{equation}
We are choosing to work with $B_\pi$ for purely expositional reasons:  it is easy to see that the transforms $\widehat{H}_p(\pi_p)$ of $1_{\chi,n}$ and $\overline{B_\pi(1_{\overline{\chi},n})}$ differ only by a multiplicative constant. Recalling \cite[Lemma 4.2]{Bessel}, this transform can be understood in terms of the action of $w$ in the Kirillov model: we have
$$B_\pi(1_{\chi,n})=\pi(w)1_{\chi_,n}(\beta)$$
\noindent In the sections to come we compute explicitly the image of $1_{\chi,n}$ under $\pi(w)$ as we vary over tuples $(\pi,\chi,n)$, providing bounds in the case of interest to us, namely $\beta \in \mathbb{Z}_p^\times$.

\subsection{Action of the Weyl element}
We use description of the action of $\pi(w)$ on the Kirillov model $\mathcal{K}(\pi,\psi)$ given in \cite[proposition 2.10]{JacquetLanglands}. Fix $f\in \mathcal{K}(\pi,\psi)$ and $\nu \in \widehat{\mathbb{Z}_p^\times}$, and let $m \in \mathbb{Z}$. We define
$$\hat{f}_m(\nu)=\int_{\mathbb{Z}_p^\times}f(u p^m)\nu(u)d^\times u \hspace{1.5cm} \hat{f}(\nu,t)=\sum_{m\in \mathbb{Z}}\hat{f}_m(\nu)t^m$$
The action of $w$ is then characterised by the Jacquet-Langlands functional equation
\begin{equation}
    \widehat{wf}(\nu,t)=C(\nu,t)\hat{f}(\nu^{-1},t^{-1})
\end{equation}
where $C(\nu,t)$ is a formal series satisfying
\begin{equation}
    C(\nu,p^{s-1/2})=\gamma(\pi,\nu^{-1},s):=\frac{L(1-s,\pi\times \nu)\varepsilon(s,\pi\times\nu^{-1},\psi)}{L(s,\pi\times\nu^{-1})}
\end{equation}
By Fourier inversion, the formal series $\widehat{wf}(\nu,t)$ determines $wf$ uniquely. Thus given $g \in \mathcal{K}(\pi,\psi)$ such that $\widehat{g}=\widehat{wf}$, it automatically implies $wf=g$.
\subsubsection{Supercuspidal representations}
An explicit expression for the action of the Weyl element is already known in the supercuspidal case (see for instance \cite[Lemma 2.1]{Saito1993}): one has
\begin{equation}\label{eq:supercuspidal-transform}
    \pi(w)1_{\chi,n} = \varepsilon(1/2,\pi\times \chi^{-1},\psi)1_{\chi^{-1},-c(\pi\times\chi^{-1})-n}
\end{equation}
We provide a similar formula for an arbitrary special representation, as well as in some cases where $\pi$ is a principal series. We treat the remaining cases using \eqref{def:induced-model-bessel} for $j_{\pi,\psi}$ to obtain bounds for $B_\pi(1_{\chi,n})$.

\subsubsection{Special representations} \label{section:bessel-special}
We now consider the representation $\pi \simeq \eta St$, where $\eta$ is a quadratic character.  Setting $f=1_{\chi,n}$, a simple computation shows that
$$\hat{f}(\nu,t)=\begin{cases}
    t^n \qquad \nu=\chi^{-1} \\
    0 \quad \hspace{0.5cm }\text{otherwise}   \end{cases}$$
thus the functional equation implies that $\widehat{wf}(\nu,t)=0$ whenever $\nu \neq \chi$. It remains to consider the case $\nu=\chi$,
we consider the following two cases:\bigbreak
\noindent 
\textbf{Case 1}: $\chi \eta \vert_{\mathbb{Z}_p^\times}$ is ramified. Since the $L$-factors are identically $1$ (using the case hypothesis) and 
$$\varepsilon(s,\pi\times \chi^{-1},\psi )=\varepsilon(1/2,\pi\times \chi^{-1},\psi)\ p^{-c(\pi\times \chi^{-1})(s-1/2)}$$
thus we obtain
$\widehat{wf}(\chi,t)=\varepsilon(1/2,\pi\times \chi^{-1},\psi)t^{-c(\pi\times \chi^{-1})-n}$. By inspection we must have
\begin{equation}
    \pi(w)1_{\chi,n} = \varepsilon(1/2,\pi\times\chi^{-1},\psi)1_{\chi^{-1},{-c(\pi\times\chi^{-1})-n}}
\end{equation} \bigbreak

\noindent \textbf{Case 2:} $\chi\eta\vert_{\mathbb{Z}_p^\times}$ is unramified. In this case
$$L(s,\pi\times \chi) = L(s,\pi\times\chi^{-1})= (1-\eta(p)p^{-s-1/2})^{-1}$$
and therefore
$$C(\chi,q^{s-1/2})=\varepsilon(1/2,\pi\times\chi^{-1},\psi)p^{-c(\pi\times \chi^{-1})(s-1/2)}\frac{1-\eta(p)p^{-(s-1/2)-1}}{1-\eta(p)p^{s-1/2-1}}$$
Noting that $c(\pi\times \chi^{-1})=1$, we may rewrite as a formal series (provided $\nu=\chi^{-1}$):
\begin{equation}
    \begin{split}
        \widehat{wf}(\nu,t ) & =C(\chi,t)\hat{f}(\chi^{-1},t^{-1}) =\varepsilon(1/2,\pi\times\chi^{-1},\psi)t^{-1-n}\frac{1-\eta(p)p^{-1}t^{-1}}{1-\eta(p)p^{-1}t} \\
        & =\varepsilon(1/2,\pi\times\chi^{-1},\psi)t^{-1-n}(1-\eta(p)p^{-1}t^{-1})\sum_{m\geq 0}\eta(p)^mp^{-m}t^{m}
    \end{split}
\end{equation}
Recall that the newform in the Kirillov model is given by $g(x)=|x|\eta(x)1_{\mathbb{Z}_p}(x)$ (see eg \cite[p.23]{newform-table}), we have
$$\widehat{g}_m(\nu)=\begin{cases}
    \eta(p)^mp^{-m} \delta_{m\geq 0}  & \text{if } \nu=\eta^{-1}\\
    0 & \text{otherwise}
\end{cases}$$
For $k \in \mathbb{Z}$, a simple check shows $\widehat{(a(p^k)g)}_m(\nu)=\widehat{g}_{m+k}(\nu)$ and therefore
\begin{equation}\label{eq:shift-whittaker}
    (\widehat{a(p^k)g})(\nu,t)=\sum_{m\geq 0}\nu(p)^m p^{-m}t^{m-k}
\end{equation}
Using \eqref{eq:shift-whittaker} taking $k=n+1$ and $k=n+2$, by linearity we deduce that
\begin{equation}
    (\pi(w)1_{\chi,n})(x) = \varepsilon(1/2,\pi\times\chi^{-1},\psi)\eta(p)^{n+1}p^{-(n+1)}|x|\left( 1_{\mathbb{Z}_p}(xp^{n+1})-p^{-2} 1_{\mathbb{Z}_p}(xp^{n+2})\right)
\end{equation} 
We collect the results above in the following lemma:
\begin{lemma}\label{lemma:bessel-function-special-bounds}
    Let $\chi$ be a character of $\mathbb{Z}_p^\times$ and $n \in \mathbb{Z}$. Suppose that $\pi\simeq \eta St$ where $\eta$ is a quadratic character and $\beta \in \mathbb{Z}_p^\times$. If $\chi\eta$ is ramified then $B_\pi(1_{\chi,n})$ vanishes unless $n=-c(\pi\times\chi^{-1})$, in which case
    $$|B_\pi(1_{\chi,n})|=1$$
    If instead $\chi\eta$ is unramified, the transform $B_\pi(1_{\chi,n})$ vanishes unless $n\geq-2$, in which case
    $$|B_\pi(1_{\chi,n})|\leq\begin{cases}
        p^{-1} & \text{if } n=-2 \\
        p^{-(n+1)} & \text{if } n>-2
    \end{cases}$$
    where the implied constants are independent of $p$.
\end{lemma}

\bigbreak 

\subsubsection{Principal series representations} Let $\mu$ be a character of $\mathbb{Q}_p^\times$, and take $\pi \simeq \pi(\mu,\mu^{-1})$. Following essentially the same argument as \S \ref{section:bessel-special}, if $\chi\mu$ and $\chi\mu^{-1}$ are both ramified, we have
\begin{equation}\label{eq:princ-series-ramified-bound}
        \pi(w)1_{\chi,n}= \varepsilon(1/2,\pi\times \chi^{-1},\psi)1_{\chi^{-1},{-c(\pi\times\chi^{-1})-n}}
\end{equation}
and the same bound as in the first part of Lemma \ref{lemma:bessel-function-special-bounds} applies.
\begin{lemma}\label{lemma:princ-series-one-ramifies-bound}
    Fix $\pi \simeq \pi(\mu,\mu^{-1})$. Let $\chi$ be a  character such that $\chi\mu \vert_{\mathbb{Z}_p^\times}$ is unramified and $\chi\mu^{-1}\vert_{\mathbb{Z}_p^\times}$ is ramified, then
    \begin{equation}
    \begin{split}
        \pi(w)1_{\chi,n}(x) & =\varepsilon(1/2,\pi\times\chi^{-1},\psi)\mu(xp^{c(\pi\times\chi^{-1})+n})p^{-(c(\pi\times\chi^{-1})+n)/2} \\
        & \quad \times |x|^{1/2}\left( 1_{\mathbb{Z}_p}(xp^{c(\pi\times\chi^{-1})+n})-p^{-1}1_{\mathbb{Z}_p}(xp^{c(\pi\times\chi^{-1})+n+1})\right)
    \end{split}
    \end{equation}
\end{lemma}

\begin{proof}
    Let $\mu_1:=\mu$ and $\mu_2:=\mu^{-1}$; we have the following $L$-factors
    $$ L(s,\pi\times\chi)=(1-\mu_1(p)p^{-s})^{-1} \quad \text{and}\quad L(s,\pi\times\chi^{-1})=(1-\mu_2(p)p^{-s})^{-1}$$
    so that
    \begin{equation}
        \begin{split}
            C(\chi,p^{s-1/2})=\varepsilon(1/2,\pi\times \chi^{-1},\psi)p^{-{c}(\pi\times\chi^{-1})(s-1/2)}\frac{1-\mu_2(p)p^{-(s-1/2)-1/2}}{1-\mu_1(p)p^{(s-1/2)-1/2}}
        \end{split}
    \end{equation}
    so given $f=1_{\chi,n}$, the functional equation implies that $\widehat{wf}(\nu,t)=0$ for $\nu \neq \chi$ and
    \begin{equation}
        \begin{split}
            \widehat{wf}(\chi,t)& =\varepsilon(1/2,\pi\times\chi^{-1},\psi)t^{-(c(\pi\times\chi^{-1})+n)}\frac{1-\mu_2(p)p^{-1/2}t^{-1}}{1-\mu_1(p)p^{-1/2}t} \\
            & =\varepsilon(1/2,\pi\times\chi^{-1},\psi)t^{-(c(\pi\times\chi^{-1})+n)}(1-\mu_2(p)p^{-1/2}t^{-1})\sum_{m\geq 0}\mu_1(p)^m p^{-m/2}t^m
        \end{split}
    \end{equation}
    Letting $g(x)=\mu_1(x)|x|^{1/2}1_{\mathbb{Z}_p}(x) \in \mathcal{K}(\pi,\psi)$ (see e.g. \cite[p.38]{Godement2018} ), a simple calculation shows that 
    $$\hat{g}_m(\chi) =\begin{cases} \mu_1(\varpi)^mp^{-m/2}\delta_{m\geq 0} \qquad \text{if } \nu=\chi\\
    0 \hspace{3.2cm} \text{otherwise}
        
    \end{cases}$$
    The result then follows after rearranging, via a computation similar to that of Lemma \ref{lemma:bessel-function-special-bounds}.
\end{proof}

\bigbreak \noindent Lastly we need to address the case $\pi\simeq \pi(\mu,\mu^{-1})$ where $\mu\vert_{\mathbb{Z}_p^\times}$ is quadratic. Although in principle the functional equation can be used to address this case as well via the power series approach, it is rather cumbersome and thus we choose instead to work with \eqref{def:induced-model-bessel}.\\
We make use of the following result:
\begin{lemma}\label{lemma:unramified-gauss-sum}
    Let $\psi$ be an additive character of conductor $\mathbb{Z}_p$ and $\chi$ an unramified multiplicative character of $\mathbb{Q}_p^\times$. Define
    $$I(m)=\int_{p^m\mathbb{Z}_p^\times}\chi^{-1}(x)\psi(x)dx $$
    then $I(m)=0$ if $m\leq-2$ and
    \begin{equation}
            I(m)=\begin{cases}
                \hspace{1cm}-\chi(p) &  \text{if} \quad m=-1 \\
                \chi(p^{-m})p^{-m}(1-p^{-1}) & \text{if} \quad m\geq 0
            \end{cases}
    \end{equation}
\end{lemma}
\begin{proof}
    This is a consequence of character orthogonality. See also \cite[Lemma 1.1.1]{newform-table}.
\end{proof}

\bigbreak \noindent We have $\pi\simeq\pi(\eta|\cdot|^{i\theta/\log p},\eta|\cdot|^{-i\theta/\log p})$ where $\eta$ is a character of $\mathbb{Q}_p^\times$ such that $\eta(p)=1$. We split the integral formula \eqref{def:induced-model-bessel} for $j_{\pi,\psi}$ into annuli of different valuation, so that
$$j_{\pi,\psi}(x)=|x|^{1/2+i\theta/\log p}\eta(x)\sum_{l \in \mathbb{Z}}J_{l}(x), \qquad J_l(x):=\int_{v(z)=l}\psi(xz+z^{-1})|z|^{2i\theta/\log p}d^\times z$$
We proceed to evaluate $J_l(x)$ for $x \in p^{-n}\mathbb{Z}_p^\times$ and $l \in \mathbb{Z}$. Suppose $n>0$, then for $l\geq n$ we have $\psi(xz)=1$; therefore if $n\geq \max(l,2)$ upon changing variables we obtain $J_l(x)=0$ using Lemma \ref{lemma:unramified-gauss-sum}.
In the range $l\leq -1$ we have instead $\psi(z^{-1})=1$ hence
$$J_l(x)=e^{-2in\theta}\int_{v(z)=l-n}\psi(z)(|z|^{-2i\theta/\log p})^{-1}d^\times z=0$$
 since by assumption we always have $l<n-1$. Thus $J_l(x)$ is supported for $0\leq l \leq n$ and applying trivial estimates
\begin{equation}\label{eq:bound-on-J_l}
    J_{l}(x)\leq \begin{cases}
        \zeta_p(1)^{-1} \qquad \text{ if } \theta \in [0,\pi] \\
        \zeta_p(1)^{-1}p^{2l\tau} \ \text{ if } \theta=i\tau\log p \text{ or } \pi+i\tau\log p \text{ where } \tau\in(0,1/2)
    \end{cases}
\end{equation}
Similarly if $n\leq 0$, the function $J_l(x)$ is supported on $n-1\leq l\leq 1$ and the same bounds as \eqref{eq:bound-on-J_l} applies. We reach the following result:
\begin{lemma}\label{lemma:princ-series-unramified-bound}
    Let $\pi \simeq \pi(\eta|\cdot|^{i\theta/\log p},\eta|\cdot|^{-i\theta/\log p})$ where $\eta$ is a quadratic character so that $\eta(p)=1$. Fix $\beta \in \mathbb{Z}_p^\times$ underlying the definition of $B_\pi$. If $\chi=\eta$ we have the bound
    $$ B_\pi(1_{\chi,n})\leq\begin{cases}
        \zeta_p(1)^{-1}(|n|+2)p^{-n/2} & \text{ if } \theta\in [0,\pi] \\
        \zeta_p(1)^{-1}(|n|+2)p^{-n/2-n\tau} & \text{ if } \theta=i\tau\log p \text{ or } \pi+i\tau\log p 
    \end{cases}$$
    where the latter case occurs only if $\eta$ is trivial.
\end{lemma}
\begin{proof}
    By the definition of $1_{\chi,n}$ we have
    \begin{align*}
        B_\pi(1_{\chi,n})& =\int_{p^n \mathbb{Z}_p^\times} 1_{\chi,n}(u)|u|^{1/2+i\theta/\log p}\eta(\beta u) \sum_{l \in \mathbb{Z}}J_l(\beta u)d^\times u
    \end{align*}
    If $\theta \in [0,\pi]$ we may bound the left hand side by
    $$B_\pi(1_{\chi,n})\leq p^{-n/2}\sum_{l\in \mathbb{Z}}\int_{p^n \mathbb{Z}_p^\times}|J_l(\beta u)|d^\times u$$
    We may restrict the $l$ summation to a finite sum and use the bounds in \eqref{eq:bound-on-J_l} to deduce the first result. Similarly, if $\theta=i\tau \log p$ or $i\tau \log p+\pi$ we have
    $$B_\pi(1_{\chi,n})\leq p^{-n/2}p^{n\tau}\sum_{l\in \mathbb{Z}}\int_{p^n \mathbb{Z}_p^\times}|J_l(\beta u)|d^\times u$$
    and using again \eqref{eq:bound-on-J_l} gives the result.
\end{proof}

\subsection{Sums of Kloosterman sums}
We will now apply the backwards formula of Theorem \ref{thm:backward-formula} to prove Corollary \ref{cor:sum-of-kloosterman}. We closely follow \cite[\S 4]{Blomer_2015}. Recall $m,n\in \mathbb{Z}$ satisfy $mn>0$ and $X\geq \sqrt{mn}$ (the Linnik range), and we let $\theta$ denote the best bound towards the generalised Ramanujan-Petersson conjecture for all places.
Given $f_\infty \in C_c^\infty(\mathbb{R}_+)$, we begin by setting
$$g_\infty(x):=\left(\frac{4\pi\sqrt{mn}}{x}\right)^{1/2}f_\infty\left( \frac{4\pi \sqrt{mn}}{xX}\right)$$
so that in the Linnik range $g_\infty(x)$ is supported on $x \asymp\Xi:=\sqrt{mn}X^{-1}\leq 1$ and satisfies the inequality of sup norms
$$||g_\infty||_\infty \asymp X^{1/2}||f_\infty||_\infty$$
We may also relate Mellin transforms of $f_\infty$ and $g_\infty$ via
$$\widehat{g}_\infty(s)=X^{1/2}\left(\frac{4\pi\sqrt{mn}}{X}\right)^s\widehat{f}_\infty(1/2-s)$$
For the transforms $\widehat{H}_\infty(\pi_\infty)$ of $g_\infty$ (recall \eqref{def:archim-backward-transform}), the power series of the Bessel function yields
\begin{equation}
    \widehat{H}_\infty(\pi_\infty)\ll ||g_\infty||_\infty \Gamma(k)^{-1}
\end{equation}
if $\pi\simeq \mathcal{D}(\kappa)$, as well as
\begin{equation}
    \widehat{H}_\infty(\pi_\infty)\ll \begin{cases}
        \frac{|\widehat{g}_\infty(2it)|+|\widehat{g}_\infty(-2it)|}{|t|^{1/2}}+\frac{|\widehat{g}_\infty(2+2it)|+|\widehat{g}_\infty(2-2it)|}{|t|^{3/2}}+\frac{||g_\infty||_\infty}{|t|^{5/2}} \quad |t|\geq 1 \\
        ||g_\infty||_\infty\Xi^{-2\theta} \quad t \in (-1,1) \cup[-i\theta,i\theta]
    \end{cases}
\end{equation}
if $\pi$ is a principal series with spectral parameter $t$
(see \cite[(4.3) and (4.4)]{Blomer_2015}). \bigbreak

\noindent For the non-archimedean test function at $p$, we make the choice $f_p=1_{1,-2k}$
 the characteristic function on the annulus of valuation $-2k$. At all other finite places $\ell\neq p$ instead we choose $f_\ell=1_{K_\ell}$ so that twisted multiplicativity of Kloosterman sums yields
$$H^S(m,n,c^2)=\prod_{\ell\neq p}H_\ell(m,n,c^2)=S(\overline{p^k}m,\overline{p^k} n, c_0)$$
where $c=c_0p^k$ and $\overline{p}$ is the multiplicative inverse of $p$ modulo $c_0$. We apply Theorem \ref{thm:backward-formula} to obtain
\begin{equation}
    \sum_{(c,p^\infty)=q}\frac{f_\infty(c/X)}{c^{1/2}}p^{2k}S(\overline{q}m,\overline{q}n,c_0)=\mathcal{M}+\mathcal{H}+\mathcal{E}
\end{equation}
where $\mathcal{H}$, $\mathcal{M}$ and $\mathcal{E}$ denote the contribution of the modular/Maass/Eisenstein forms in Theorem \ref{thm:backward-formula} (with transforms with respect to $g_\infty$ and $f_p$). Corollary \ref{cor:sum-of-kloosterman} follows by bounding the spectral terms individually. For the cuspidal spectrum we limit ourselves to the Maass case, since the modular forms case follows similarly. Combining equations \eqref{eq:supercuspidal-transform}, \eqref{eq:princ-series-ramified-bound} and Lemmata \ref{lemma:bessel-function-special-bounds}, \ref{lemma:princ-series-one-ramifies-bound} and \ref{lemma:princ-series-unramified-bound}, for $\pi_p \in G(\mathbb{Q}_p)^\wedge$ and $k$ sufficiently large, we have $M_{\pi_p}(f_p)=0$ unless $c(\pi_p)=2k$, in which case
\begin{equation}\label{eq:high-conductor}
    \widehat{H}_{p}(\pi_p) \ll 1,
\end{equation}
or $c(\pi_p)=0$, for which we may obtain
\begin{equation}\label{eq:unram-rep}
    \widehat{H}_{p}(\pi_p) \ll p^{k(1+2\tau)}\log p^k.
\end{equation}
where $\tau$ is the largest bound in the exceptional spectrum of $G(\mathbb{Q}_q)^\wedge$ appearing as a component of a cuspidal automorphic representation. As mentioned in \cite[(1.12)]{generalised}, we have  that for a normalised newform $\varphi \in V_\pi$:  $$a_{u_\varphi}(m)\overline{a_{u_\varphi}(n)}\ll_\varepsilon (C(\pi_p)(1+|t_\pi|))^\epsilon (|mn|)^{\theta+\epsilon}$$
By the Weyl's law for Maass forms on $\Gamma_0(p^{2k})$ and the bounds \eqref{eq:high-conductor} and \eqref{eq:unram-rep}, we obtain
\begin{equation}\label{eq:maass-forms-bound}
    \mathcal{M}\ll_{f_\infty,\epsilon} X^{1/2+2\theta}p^{(2+\epsilon)k}(mn)^\epsilon
\end{equation}
where, explicitly, the dependence on $f_\infty$ is given by
\begin{align}\label{eq:f-infinity-dependence}
    \# \{\pi \ \text{spherical}: & |t_\pi|<1 \} \cdot ||f_\infty||_\infty \notag \\ 
    &+\sum _{\pi \ \text{spherical}}\left( \frac{|\widehat{f}_\infty(1/2\pm2it_\pi)|}{|t_\pi|^{1/2-\epsilon}}+\frac{|\widehat{f}_\infty(-3/2\pm2it_\pi)}{|t_\pi|^{3/2-\epsilon}}+\frac{||f_\infty||_\infty}{|t_\pi|^{5/2-\epsilon}}\right)
\end{align}
with the summation being over all cuspidal spherical representations. As mentioned previously, the holomorphic case is similar, we obtain
\begin{equation*}
    \mathcal{H}\ll_\epsilon ||f_\infty||_\infty X^{1/2}p^{2(k+\epsilon)}(mn)^\epsilon 
\end{equation*}
Lastly in the Eisenstein series case, we have
$$ a_{u_{E(\phi_{it})}}(m)\overline{a_{u_{E(\phi_{it})}}(n)}\ll_\epsilon (C(\pi_p)(1+|t|))^\epsilon (mn)^\epsilon $$
which holds for all induced model representations. We obtain similar bounds as in \eqref{eq:maass-forms-bound}, replacing \eqref{eq:f-infinity-dependence} with
\begin{equation*}
    ||f_\infty||_\infty +\int_{1}^{\infty}\frac{|\widehat{f}_\infty(1/2-2it)|}{|t|^{1/2-\epsilon}}dt
\end{equation*} 
Collecting all terms then proves Corollary \ref{cor:sum-of-kloosterman}.

\begin{proof}[Proof of Corollary \ref{cor:dedekind}]
    From \cite[Theorem 1.1]{vardi-1987} we have
    \begin{equation}
        \sum_{\substack{0<d<c \\(d,c)=1}}e(12ms(d,c))=S(m,m,c)
    \end{equation}
    therefore we have from Corollary \ref{cor:sum-of-kloosterman},
    \begin{align}\label{eq:non-trivial-dedekind-bound}
        \sum_{\substack{c\in \mathbb{N}\\ (c,p^\infty)=q}}f_\infty(c/X)\sum_{\substack{0<d<c \\ (c,d)=1}}e(12\overline{q}s(d,c_0))\ll_{f_\infty} X^{1+2\theta+\epsilon}q^{\epsilon}
    \end{align}
    Since the bound in \eqref{eq:non-trivial-dedekind-bound} is not trivial, the equidistribution results follows from Weyl equiditribution criterion upon replacing $f_\infty$ with a characteristic function on $[0,1]$ (eg via diadic partition).
\end{proof}

\pagebreak

\pagebreak 
\printbibliography
\end{document}